# KERNEL DIMENSION REDUCTION IN REGRESSION[1]

By Kenji Fukumizu, Francis R. Bach and Michael I. Jordan

*Institute of Statistical Mathematics, INRIA—Ecole Normale Supérieure and University of California*

We present a new methodology for sufficient dimension reduction (SDR). Our methodology derives directly from the formulation of SDR in terms of the conditional independence of the covariate $X$ from the response $Y$, given the projection of $X$ on the central subspace [cf. *J. Amer. Statist. Assoc.* **86** (1991) 316–342 and *Regression Graphics* (1998) Wiley]. We show that this conditional independence assertion can be characterized in terms of conditional covariance operators on reproducing kernel Hilbert spaces and we show how this characterization leads to an $M$-estimator for the central subspace. The resulting estimator is shown to be consistent under weak conditions; in particular, we do not have to impose linearity or ellipticity conditions of the kinds that are generally invoked for SDR methods. We also present empirical results showing that the new methodology is competitive in practice.

**1. Introduction.** The problem of *sufficient dimension reduction* (SDR) for regression is that of finding a subspace $S$ such that the projection of the covariate vector $X$ onto $S$ captures the statistical dependency of the response $Y$ on $X$. More formally, let us characterize a *dimension-reduction subspace* $S$ in terms of the following conditional independence assertion:

$$(1) \qquad Y \perp\!\!\!\perp X | \Pi_S X,$$

where $\Pi_S X$ denotes the orthogonal projection of $X$ onto $S$. It is possible to show that under weak conditions the intersection of dimension-reduction subspaces is itself a dimension-reduction subspace, in which case the intersection is referred to as a *central subspace* [5, 6]. As suggested in a seminal

Received October 2006; revised July 2008.

[1]Supported by JSPS KAKENHI 15700241, a grant from the Inamori Foundation, a Scientific grant from the Mitsubishi Foundation and NSF Grant DMS-05-09559.

paper by Li [23], it is of great interest to develop procedures for estimating this subspace, quite apart from any interest in the conditional distribution $P(Y|X)$ or the conditional mean $E(Y|X)$. Once the central subspace is identified, subsequent analysis can attempt to infer a conditional distribution or a regression function using the (low-dimensional) coordinates $\Pi_S X$.

The line of research on SDR initiated by Li is to be distinguished from the large and heterogeneous collection of methods for dimension reduction in regression in which specific modeling assumptions are imposed on the conditional distribution $P(Y|X)$ or the regression $E(Y|X)$. These methods include ordinary least squares, partial least squares, canonical correlation analysis, ACE [4], projection pursuit regression [12], neural networks and LASSO [29]. These methods can be effective if the modeling assumptions that they embody are met, but if these assumptions do not hold there is no guarantee of finding the central subspace.

Li's paper not only provided a formulation of SDR as a semiparametric inference problem—with subsequent contributions by Cook and others bringing it to its elegant expression in terms of conditional independence—but also suggested a specific inferential methodology that has had significant influence on the ensuing literature. Specifically, Li suggested approaching the SDR problem as an *inverse regression* problem. Roughly speaking, the idea is that if the conditional distribution $P(Y|X)$ varies solely along a subspace of the covariate space, then the inverse regression $E(X|Y)$ should lie in that same subspace. Moreover, it should be easier to regress $X$ on $Y$ than vice versa, given that $Y$ is generally low-dimensional (indeed, one-dimensional in the majority of applications) while $X$ is high-dimensional. Li [23] proposed a particularly simple instantiation of this idea—known as *sliced inverse regression* (SIR)—in which $E(X|Y)$ is estimated as a constant vector within each slice of the response variable $Y$, and principal component analysis is used to aggregate these constant vectors into an estimate of the central subspace. The past decade has seen a number of further developments in this vein. Some focus on finding a central subspace, for example, [9, 10], while others aim at finding a *central mean subspace*, which is a subspace of the central subspace that is effective only for the regression $E[Y|X]$. The latter include principal Hessian directions (pHd, [24]) and contour regression [22]. A particular focus of these more recent developments has been the exploitation of second moments within an inverse regression framework.

While the inverse regression perspective has been quite useful, it is not without its drawbacks. In particular, performing a regression of $X$ on $Y$ generally requires making assumptions with respect to the probability distribution of $X$, assumptions that can be difficult to justify. In particular, most of the inverse regression methods make the assumption of linearity of the conditional mean of the covariate along the central subspace (or make a related assumption for the conditional covariance). These assumptions hold in particular if the distribution of $X$ is elliptic. In practice, however, we do not necessarily expect that the covariate vector will follow an elliptic distribution, nor is it easy to assess departures from ellipticity in a high-dimensional setting. In general, it seems unfortunate to have to impose probabilistic assumptions on $X$ in the setting of a regression methodology.

Many of inverse regression methods can also exhibit some additional limitations depending on the specific nature of the response variable $Y$. In particular, pHd and contour regression are applicable only to a one-dimensional response. Also, if the response variable takes its values in a finite set of $p$ elements, SIR yields a subspace of dimension at most $p - 1$; thus, for the important problem of binary classification SIR yields only a one-dimensional subspace. Finally, in the binary classification setting, if the covariance matrices of the two classes are the same, SAVE and pHd also provide only a one-dimensional subspace [7]. The general problem in these cases is that the estimated subspace is smaller than the central subspace. One approach to tackling these limitations is to incorporate higher-order moments of $Y|X$ [34], but in practice the gains achievable by the use of higher-order moments are limited by robustness issues.

In this paper, we present a new methodology for SDR that is rather different from the approaches considered in the literature discussed above. Rather than focusing on a limited set of moments within an inverse regression framework, we focus instead on the criterion of conditional independence in terms of which the SDR problem is defined. We develop a contrast function for evaluating subspaces that is minimized precisely when the conditional independence assertion in (1) is realized. As befits a criterion that measures departure from conditional independence, our contrast function is not based solely on low-order moments.

Our approach involves the use of conditional covariance operators on reproducing kernel Hilbert spaces (RKHSs). Our use of RKHSs is related to their use in nonparametric regression and classification; in particular, the RKHSs given by some positive definite kernels are Hilbert spaces of smooth functions that are "small" enough to yield computationally-tractable procedures, but are rich enough to capture nonparametric phenomena of interest [32], and this computational focus is an important aspect of our work. On the other hand, whereas in nonparametric regression and classification the role of RKHSs is to provide basis expansions of regression functions and discriminant functions, in our case the RKHS plays a different role. Our interest is not in the functions in the RKHS per se, but rather in conditional covariance operators defined on the RKHS. We show that these operators can be used to measure departures from conditional independence. We also show that these operators can be estimated from data and that these estimates are functions of Gram matrices. Thus, our approach—which we refer to as *kernel dimension reduction* (KDR)—involves computing Gram matrices from data and optimizing a particular functional of these Gram matrices to yield an estimate of the central subspace.

This approach makes no strong assumptions on either the conditional distribution $p_{Y|\Pi_S X}(y|\Pi_S x)$ or the marginal distribution $p_X(x)$. As we show, KDR is consistent as an estimator of the central subspace under weak conditions.

There are alternatives to the inverse regression approach in the literature that have some similarities to KDR. In particular, minimum average variance estimation (MAVE, [33]) is based on nonparametric estimation of the conditional covariance of $Y$ given $X$, an idea related to KDR. This method explicitly estimates the regressor, however, assuming an additive noise model $Y = f(X) + Z$, where $Z$ is independent of $X$. While the purpose of MAVE is to find a central mean subspace, KDR tries to find a central subspace, and does not need to estimate the regressor explicitly. Other related approaches include methods that estimate the derivative of the regression function; these are based on the fact that the derivative of the conditional expectation $g(x) = E[y|B^T x]$ with respect to $x$ belongs to a dimension reduction subspace [18, 27]. The purpose of these methods is again to extract a central mean subspace; this differs from the central subspace which is the focus of KDR. The difference is clear, for example, if we consider the situation in which a direction $b$ in a central subspace satisfies $E[g'(b^T X)] = 0$; a condition that occurs if $g$ and the distribution of $X$ exhibit certain symmetries. The direction cannot be found by methods based on the derivative. Also, there has also been some recent work on nonparametric methods for estimation of central subspaces. One such method estimates the central subspace based on an expected log likelihood [35]. This requires, however, an estimate of the joint probability density, and is limited to single-index regression. Finally, Zhu and Zeng [36] have proposed a method for estimating the central subspace based on the Fourier transform. This method is similar to the KDR method in its use of Hilbert space methods and in its use of a contrast function that can characterize conditional independence under weak assumptions. It differs from KDR, however, in that it requires an estimate of the derivative of the marginal density of the covariate $X$; in practice this requires assuming a parametric model for the covariate $X$. In general, we are aware of no practical method that attacks SDR directly by using nonparametric methodology to assess departures from conditional independence.

We presented an earlier kernel dimension reduction method in [14]. The contrast function presented in that paper, however, was not derived as an estimator of a conditional covariance operator, and it was not possible to establish a consistency result for that approach. The contrast function that we present here is derived directly from the conditional covariance perspective;

moreover, it is simpler than the earlier estimator and it is possible to establish consistency for the new formulation. We should note, however, that the empirical performance of the earlier KDR method was shown by Fukumizu, Bach and Jordan [14] to yield a significant improvement on SIR and pHd in the case of nonelliptic data, and these empirical results motivated us to pursue the general approach further.

While KDR has advantages over other SDR methods because of its generality and its directness in capturing the semiparametric nature of the SDR problem, it also reposes on a more complex mathematical framework that presents new theoretical challenges. Thus, while consistency for SIR and related methods follows from a straightforward appeal to the central limit theorem (under ellipticity assumptions), more effort is required to study the statistical behavior of KDR theoretically. This effort is of some general value, however; in particular, to establish the consistency of KDR we prove the uniform $O(n^{-1/2})$ convergence of an empirical process that takes values in a reproducing kernel Hilbert space. This result, which accords with the order of uniform convergence of an ordinary real-valued empirical process, may be of independent theoretical interest.

It should be noted at the outset that we do not attempt to provide distribution theory for KDR in this paper, and in particular we do not address the problem of inferring the dimensionality of the central subspace.

The paper is organized as follows. In Section 2 we show how conditional independence can be characterized by cross-covariance operators on an RKHS and use this characterization to derive the KDR method. Section 3 presents numerical examples of the KDR method. We present a consistency theorem and its proof in Section 4. Section 5 provides concluding remarks. Some of the details in the proof of consistency are provided in the Appendix.

**2. Kernel dimension reduction for regression.** The method of kernel dimension reduction is based on a characterization of conditional independence using operators on RKHSs. We present this characterization in Section 2.1 and show how it yields a population criterion for SDR in Section 2.2. This population criterion is then turned into a finite-sample estimation procedure in Section 2.3.

In this paper, a Hilbert space means a separable Hilbert space, and an operator always means a linear operator. The operator norm of a bounded operator $T$ is denoted by $\|T\|$. The null space and the range of an operator $T$ are denoted by $\mathcal{N}(T)$ and $\mathcal{R}(T)$, respectively.

2.1. *Characterization of conditional independence.* Let $(\mathcal{X}, \mathcal{B}_\mathcal{X})$ and $(\mathcal{Y}, \mathcal{B}_\mathcal{Y})$ denote measurable spaces. When the base space is a topological space, the Borel $\sigma$-field is always assumed. Let $(\mathcal{H}_\mathcal{X}, k_\mathcal{X})$ and $(\mathcal{H}_\mathcal{Y}, k_\mathcal{Y})$ be RKHSs of functions on $\mathcal{X}$ and $\mathcal{Y}$, respectively, with measurable positive definite kernels $k_\mathcal{X}$ and $k_\mathcal{Y}$ [1]. We consider a random vector $(X,Y):\Omega \to \mathcal{X} \times \mathcal{Y}$ with the law $P_{XY}$. The marginal distribution of $X$ and $Y$ are denoted by $P_X$ and $P_Y$, respectively. It is always assumed that the positive definite kernels satisfy

$$(2) \qquad E_X[k_\mathcal{X}(X,X)] < \infty \quad \text{and} \quad E_Y[k_\mathcal{Y}(Y,Y)] < \infty.$$

Note that any bounded kernels satisfy this assumption. Also, under this assumption, $\mathcal{H}_\mathcal{X}$ and $\mathcal{H}_\mathcal{Y}$ are included in $L^2(P_X)$ and $L^2(P_Y)$, respectively, where $L^2(\mu)$ denotes the Hilbert space of square integrable functions with respect to the measure $\mu$, and the inclusions $J_\mathcal{X}:\mathcal{H}_\mathcal{X} \to L^2(P_X)$ and $J_\mathcal{Y}:\mathcal{H}_\mathcal{Y} \to L^2(P_Y)$ are continuous, because $E_X[f(X)^2] = E_X[\langle f, k_\mathcal{X}(\cdot, X)\rangle^2_{\mathcal{H}_\mathcal{X}}] \leq \|f\|^2_{\mathcal{H}_\mathcal{X}} E_X[k_\mathcal{X}(X,X)]$ for $f \in \mathcal{H}_\mathcal{X}$.

The *cross-covariance operator* of $(X,Y)$ is an operator from $\mathcal{H}_\mathcal{X}$ to $\mathcal{H}_\mathcal{Y}$ so that

$$(3) \qquad \langle g, \Sigma_{YX} f\rangle_{\mathcal{H}_\mathcal{Y}} = E_{XY}[(f(X) - E_X[f(X)])(g(Y) - E_Y[g(Y)])]$$

holds for all $f \in \mathcal{H}_\mathcal{X}$ and $g \in \mathcal{H}_\mathcal{Y}$ [3, 14]. Obviously, $\Sigma_{YX} = \Sigma^*_{XY}$, where $T^*$ denotes the adjoint of an operator $T$. If $Y$ is equal to $X$, the positive self-adjoint operator $\Sigma_{XX}$ is called the *covariance operator*.

For a random variable $X:\Omega \to \mathcal{X}$, the *mean element* $m_X \in \mathcal{H}_\mathcal{X}$ is defined by the element that satisfies

$$(4) \qquad \langle f, m_X \rangle_{\mathcal{H}_\mathcal{X}} = E_X[f(X)]$$

for all $f \in \mathcal{H}_\mathcal{X}$; that is, $m_X = J^*_\mathcal{X} 1$, where 1 is the constant function. The explicit function form of $m_X$ is given by $m_X(u) = \langle m_X, k(\cdot, u)\rangle_{\mathcal{H}_\mathcal{X}} = E[k(X,u)]$. Using the mean elements, (3), which characterizes $\Sigma_{YX}$, can be written as

$$\langle g, \Sigma_{YX} f\rangle_{\mathcal{H}_\mathcal{Y}} = E_{XY}[\langle f, k_\mathcal{X}(\cdot, X) - m_X\rangle_{\mathcal{H}_\mathcal{X}} \langle k_\mathcal{Y}(\cdot, Y) - m_Y, g\rangle_{\mathcal{H}_\mathcal{Y}}].$$

Let $Q_X$ and $Q_Y$ be the orthogonal projections which map $\mathcal{H}_\mathcal{X}$ onto $\overline{\mathcal{R}(\Sigma_{XX})}$ and $\mathcal{H}_\mathcal{Y}$ onto $\overline{\mathcal{R}(\Sigma_{YY})}$, respectively. It is known [3], Theorem 1, that $\Sigma_{YX}$ has a representation of the form

$$(5) \qquad \Sigma_{YX} = \Sigma^{1/2}_{YY} V_{YX} \Sigma^{1/2}_{XX},$$

where $V_{YX}:\mathcal{H}_\mathcal{X} \to \mathcal{H}_\mathcal{Y}$ is a unique bounded operator such that $\|V_{YX}\| \leq 1$ and $V_{YX} = Q_Y V_{YX} Q_X$.

A cross-covariance operator on an RKHS can be represented explicitly as an integral operator. For arbitrary $\varphi \in L^2(P_X)$ and $y \in \mathcal{Y}$, the integral

$$(6) \qquad G_\varphi(y) = \int_{\mathcal{X} \times \mathcal{Y}} k_\mathcal{Y}(y, \tilde{y})(\varphi(\tilde{x}) - E_X[\varphi(X)]) \, dP_{XY}(\tilde{x}, \tilde{y})$$

always exists and $G_\varphi$ is an element of $L^2(P_Y)$. It is not difficult to see that

$$S_{YX}:L^2(P_X) \to L^2(P_Y), \qquad \varphi \mapsto G_\varphi$$

is a bounded linear operator with $\|S_{YX}\| \leq E_Y[k_{\mathcal{Y}}(Y,Y)]$. If $f$ is a function in $\mathcal{H}_{\mathcal{X}}$, we have for any $y \in \mathcal{Y}$

$$G_f(y) = \langle k_{\mathcal{Y}}(\cdot, y), \Sigma_{YX} f \rangle_{\mathcal{H}_{\mathcal{Y}}} = (\Sigma_{YX} f)(y),$$

which implies the following proposition:

PROPOSITION 1. *The covariance operator $\Sigma_{YX}: \mathcal{H}_{\mathcal{X}} \to \mathcal{H}_{\mathcal{Y}}$ is the restriction of the integral operator $S_{YX}$ to $\mathcal{H}_{\mathcal{X}}$. More precisely,*

$$J_{\mathcal{Y}} \Sigma_{YX} = S_{YX} J_{\mathcal{X}}.$$

Conditional variance can be also represented by covariance operators. Define the *conditional covariance operator* $\Sigma_{YY|X}$ by

$$\Sigma_{YY|X} = \Sigma_{YY} - \Sigma_{YY}^{1/2} V_{YX} V_{XY} \Sigma_{YY}^{1/2},$$

where $V_{YX}$ is the bounded operator in (5). For convenience we sometimes write $\Sigma_{YY|X}$ as

$$\Sigma_{YY|X} = \Sigma_{YY} - \Sigma_{YX} \Sigma_{XX}^{-1} \Sigma_{XY},$$

which is an abuse of notation, because $\Sigma_{XX}^{-1}$ may not exist.

The following two propositions provide insights into the meaning of a conditional covariance operator. The former proposition relates the operator to the residual error of regression, and the latter proposition expresses the residual error in terms of the conditional variance.

PROPOSITION 2. *For any $g \in \mathcal{H}_{\mathcal{Y}}$,*

$$\langle g, \Sigma_{YY|X} g \rangle_{\mathcal{H}_{\mathcal{Y}}} = \inf_{f \in \mathcal{H}_{\mathcal{X}}} E_{XY} |(g(Y) - E_Y[g(Y)]) - (f(X) - E_X[f(X)])|^2.$$

PROOF. Let $\Sigma_{YX} = \Sigma_{YY}^{1/2} V_{YX} \Sigma_{XX}^{1/2}$ be the decomposition in (5), and define $\mathcal{E}_g(f) = E_{YX}|(g(Y) - E_Y[g(Y)]) - (f(X) - E_X[f(X)])|^2$. From the equality

$$\mathcal{E}_g(f) = \|\Sigma_{XX}^{1/2} f\|_{\mathcal{H}_{\mathcal{X}}}^2 - 2\langle V_{XY} \Sigma_{YY}^{1/2} g, \Sigma_{XX}^{1/2} f \rangle_{\mathcal{H}_{\mathcal{X}}} + \|\Sigma_{YY}^{1/2} g\|_{\mathcal{H}_{\mathcal{Y}}}^2,$$

replacing $\Sigma_{XX}^{1/2} f$ with an arbitrary $\phi \in \mathcal{H}_{\mathcal{X}}$ yields

$$\inf_{f \in \mathcal{H}_{\mathcal{X}}} \mathcal{E}_g(f) \geq \inf_{\phi \in \mathcal{H}_{\mathcal{X}}} \{\|\phi\|_{\mathcal{H}_{\mathcal{X}}}^2 - 2\langle V_{XY} \Sigma_{YY}^{1/2} g, \phi \rangle_{\mathcal{H}_{\mathcal{X}}} + \|\Sigma_{YY}^{1/2} g\|_{\mathcal{H}_{\mathcal{Y}}}^2\}$$

$$= \inf_{\phi \in \mathcal{H}_{\mathcal{X}}} \|\phi - V_{XY} \Sigma_{YY}^{1/2} g\|_{\mathcal{H}_{\mathcal{X}}}^2 + \langle g, \Sigma_{YY|X} g \rangle_{\mathcal{H}_{\mathcal{Y}}}$$

$$= \langle g, \Sigma_{YY|X} g \rangle_{\mathcal{H}_{\mathcal{Y}}}.$$

For the opposite inequality, take an arbitrary $\varepsilon > 0$. From the fact that $V_{XY}\Sigma_{YY}^{1/2}g \in \overline{\mathcal{R}(\Sigma_{XX})} = \mathcal{R}(\Sigma_{XX}^{1/2})$, there exists $f_* \in \mathcal{H}_\mathcal{X}$ such that $\|\Sigma_{XX}^{1/2}f_* - V_{XY}\Sigma_{YY}^{1/2}g\|_{\mathcal{H}_\mathcal{X}} \leq \varepsilon$. For such $f_*$,

$$\begin{aligned}\mathcal{E}_g(f_*) &= \|\Sigma_{XX}^{1/2}f_*\|_{\mathcal{H}_\mathcal{X}}^2 - 2\langle V_{XY}\Sigma_{YY}^{1/2}g, \Sigma_{XX}^{1/2}f_*\rangle_{\mathcal{H}_\mathcal{X}} + \|\Sigma_{YY}^{1/2}g\|_{\mathcal{H}_\mathcal{Y}}^2 \\ &= \|\Sigma_{XX}^{1/2}f_* - V_{YX}\Sigma_{YY}^{1/2}g\|_{\mathcal{H}_\mathcal{X}}^2 + \|\Sigma_{YY}^{1/2}g\|_{\mathcal{H}_\mathcal{Y}}^2 - \|V_{XY}\Sigma_{YY}^{1/2}g\|_{\mathcal{H}_\mathcal{X}}^2 \\ &\leq \langle g, \Sigma_{YY|X}g\rangle_{\mathcal{H}_\mathcal{Y}} + \varepsilon^2.\end{aligned}$$

Because $\varepsilon$ is arbitrary, we have $\inf_{f \in \mathcal{H}_\mathcal{X}} \mathcal{E}_g(f) \leq \langle g, \Sigma_{YY|X}g\rangle_{\mathcal{H}_\mathcal{Y}}$. □

Proposition 2 is an analog for operators of a well-known result on covariance matrices and linear regression: the conditional covariance matrix $C_{YY|X} = C_{YY} - C_{YX}C_{XX}^{-1}C_{XY}$ expresses the residual error of the least square regression problem as $b^T C_{YY|X} b = \min_a E\|b^T Y - a^T X\|^2$.

To relate the residual error in Proposition 2 to the conditional variance of $g(Y)$ given $X$, we make the following mild assumption:

(AS) $\mathcal{H}_\mathcal{X} + \mathbb{R}$ is dense in $L^2(P_X)$, where $\mathcal{H}_\mathcal{X} + \mathbb{R}$ denotes the direct sum of the RKHS $\mathcal{H}_\mathcal{X}$ and the RKHS $\mathbb{R}$ [1].

As seen later in Section 2.2, there are many positive definite kernels that satisfy the assumption (AS). Examples include the Gaussian radial basis function (RBF) kernel $k(x,y) = \exp(-\|x-y\|^2/\sigma^2)$ on $\mathbb{R}^m$ or on a compact subset of $\mathbb{R}^m$.

PROPOSITION 3. *Under the assumption (AS),*

(7) $$\langle g, \Sigma_{YY|X}g\rangle_{\mathcal{H}_\mathcal{Y}} = E_X[\text{Var}_{Y|X}[g(Y)|X]]$$

*for all $g \in \mathcal{H}_\mathcal{Y}$.*

PROOF. From Proposition 2, we have

$$\begin{aligned}\langle g, &\Sigma_{YY|X}g\rangle_{\mathcal{H}_\mathcal{Y}} \\ &= \inf_{f \in \mathcal{H}_\mathcal{X}} \text{Var}[g(Y) - f(X)] \\ &= \inf_{f \in \mathcal{H}_\mathcal{X}} \{\text{Var}_X[E_{Y|X}[g(Y) - f(X)|X]] + E_X[\text{Var}_{Y|X}[g(Y) - f(X)|X]]\} \\ &= \inf_{f \in \mathcal{H}_\mathcal{X}} \text{Var}_X[E_{Y|X}[g(Y)|X] - f(X)] + E_X[\text{Var}_{Y|X}[g(Y)|X]].\end{aligned}$$

Let $\varphi(x) = E_{Y|X}[g(Y)|X=x]$. Since $\varphi \in L^2(P_X)$ from $\text{Var}[\varphi(X)] \leq \text{Var}[g(Y)] < \infty$, the assumption (AS) implies that for an arbitrary $\varepsilon > 0$ there exists $f \in \mathcal{H}_\mathcal{X}$ and $c \in \mathbb{R}$ such that $h = f + c$ satisfies $\|\varphi - h\|_{L^2(P_X)} < \varepsilon$. Because $\text{Var}[\varphi(X) - f(X)] \leq \|\varphi - h\|_{L^2(P_X)}^2 \leq \varepsilon^2$ and $\varepsilon$ is arbitrary, we have $\inf_{f \in \mathcal{H}_\mathcal{X}} \text{Var}_X[E_{Y|X}[g(Y)|X] - f(X)] = 0$, which completes the proof. □

Proposition 3 improves a result due to Fukumizu, Bach and Jordan [14], Proposition 5, where the much stronger assumption $E[g(Y)|X = \cdot] \in \mathcal{H}_\mathcal{X}$ was imposed.

Propositions 2 and 3 imply that the operator $\Sigma_{YY|X}$ can be interpreted as capturing the predictive ability for $Y$ of the explanatory variable $X$.

2.2. *Criterion of kernel dimension reduction.* Let $M(m \times n; \mathbb{R})$ be the set of real-valued $m \times n$ matrices. For a natural number $d \leq m$, the Stiefel manifold $\mathbb{S}_d^m(\mathbb{R})$ is defined by

$$\mathbb{S}_d^m(\mathbb{R}) = \{B \in M(m \times d; \mathbb{R}) | B^T B = I_d\},$$

which is the set of all $d$ orthonormal vectors in $\mathbb{R}^m$. It is well known that $\mathbb{S}_d^m(\mathbb{R})$ is a compact smooth manifold. For $B \in \mathbb{S}_d^m(\mathbb{R})$, the matrix $BB^T$ defines an orthogonal projection of $\mathbb{R}^m$ onto the $d$-dimensional subspace spanned by the column vectors of $B$. Although the Grassmann manifold is often used in the study of sets of subspaces in $\mathbb{R}^m$, we find the Stiefel manifold more convenient as it allows us to use matrix notation explicitly.

Hereafter, $\mathcal{X}$ is assumed to be either a closed ball $D_m(r) = \{x \in \mathbb{R}^m | \|x\| \leq r\}$ or the entire Euclidean space $\mathbb{R}^m$; both assumptions satisfy the condition that the projection $BB^T \mathcal{X}$ is included in $\mathcal{X}$ for all $B \in \mathbb{S}_d^m(\mathbb{R})$.

Let $\mathbb{B}_d^m \subseteq \mathbb{S}_d^m(\mathbb{R})$ denote the subset of matrices whose columns span a dimension-reduction subspace; for each $B_0 \in \mathbb{B}_d^m$, we have

(8)  $$p_{Y|X}(y|x) = p_{Y|B_0^T X}(y|B_0^T x),$$

where $p_{Y|X}(y|x)$ and $p_{Y|B^T X}(y|u)$ are the conditional probability densities of $Y$ given $X$, and $Y$ given $B^T X$, respectively. The existence and positivity of these conditional probability densities are always assumed hereafter. As we have discussed in the Introduction, under conditions given by [6], Section 6.4, this subset represents the central subspace (under the assumption that $d$ is the minimum dimensionality of the dimension reduction subspaces).

We now turn to the key problem of characterizing the subset $\mathbb{B}_d^m$ using conditional covariance operators on reproducing kernel Hilbert spaces. In the following, we assume that $k_d(z, \tilde{z})$ is a positive definite kernel on $\mathcal{Z} = D_d(r)$ or $\mathbb{R}^d$ such that $E_X[k_d(B^T X, B^T X)] < \infty$ for all $B \in \mathbb{S}_d^m(\mathbb{R})$, and we let $k_\mathcal{X}^B$ denote a positive definite kernel on $\mathcal{X}$ given by

(9)  $$k_\mathcal{X}^B(x, \tilde{x}) = k_d(B^T x, B^T \tilde{x})$$

for each $B \in \mathbb{S}_d^m(\mathbb{R})$. The RKHS associated with $k_{\mathcal{X}}^B$ is denoted by $\mathcal{H}_{\mathcal{X}}^B$. Note that $\mathcal{H}_{\mathcal{X}}^B = \{f : \mathcal{X} \to \mathbb{R} | \text{there exists } g \in \mathcal{H}_{k_d} \text{ such that } f(x) = g(B^T x)\}$, where $\mathcal{H}_{k_d}$ is the RKHS given by $k_d$. As seen later in Theorem 4, if $\mathcal{X}$ and $\mathcal{Y}$ are subsets of Euclidean spaces and Gaussian RBF kernels are used for $k_{\mathcal{X}}$ and $k_{\mathcal{Y}}$, under some conditions the subset $\mathbb{B}_d^m$ is characterized by the set of solutions of an optimization problem

$$\mathbb{B}_d^m = \underset{B \in \mathbb{S}_d^m(\mathbb{R})}{\arg\min} \Sigma_{YY|X}^B, \tag{10}$$

where $\Sigma_{YX}^B$ and $\Sigma_{XX}^B$ denote the (cross-) covariance operators with respect to the kernel $k^B$, and

$$\Sigma_{YY|X}^B = \Sigma_{YY} - \Sigma_{YX}^B \Sigma_{XX}^{B}{}^{-1} \Sigma_{XY}^B.$$

The minimization in (10) refers to the minimal operators in the partial order of self-adjoint operators.

We use the trace to evaluate the partial order of self-adjoint operators. While other possibilities exist (e.g., the determinant), the trace has the advantage of yielding a relatively simple theoretical analysis, which is conducted in Section 4. The operator $\Sigma_{YY|X}^B$ is trace class for all $B \in \mathbb{S}_d^m(\mathbb{R})$, since $\Sigma_{YY|X}^B \leq \Sigma_{YY}$ and $\text{Tr}[\Sigma_{YY}] < \infty$, which is shown in Section 4.2. Henceforth the minimization in (10) should thus be understood as that of minimizing $\text{Tr}[\Sigma_{YY|X}^B]$.

From Propositions 2 and 3, minimization of $\text{Tr}[\Sigma_{YY|X}^B]$ is equivalent to the minimization of the sum of the residual errors for the optimal prediction of functions of $Y$ using $B^T X$, where the sum is taken over a complete orthonormal system $\{\xi_a\}_{a=1}^{\infty}$ of $\mathcal{H}_{\mathcal{Y}}$. Thus, the objective of dimension reduction is rewritten as

$$\min_{B \in \mathbb{S}_d^m(\mathbb{R})} \sum_{a=1}^{\infty} \min_{f \in \mathcal{H}_{\mathcal{X}}^B} E|(\xi_a(Y) - E[\xi_a(Y)]) - (f(X) - E[f(X)])|^2. \tag{11}$$

This is intuitively reasonable as a criterion of choosing $B$, and we will see that this is equivalent to finding the central subspace under some conditions.

We now introduce a class of kernels to characterize conditional independence. Let $(\Omega, \mathcal{B})$ be a measurable space, let $(\mathcal{H}, k)$ be an RKHS over $\Omega$ with the kernel $k$ measurable and bounded, and let $\mathcal{S}$ be the set of all probability measures on $(\Omega, \mathcal{B})$. The RKHS $\mathcal{H}$ is called *characteristic* (with respect to $\mathcal{B}$) if the map

$$\mathcal{S} \ni P \mapsto m_P = E_{X \sim P}[k(\cdot, X)] \in \mathcal{H} \tag{12}$$

is one-to-one, where $m_P$ is the mean element of the random variable with law $P$. It is easy to see that $\mathcal{H}$ is characteristic if and only if the equality $\int f\,dP = \int f\,dQ$ for all $f \in \mathcal{H}$ means $P = Q$. We also call a positive definite kernel $k$ characteristic if the associated RKHS is characteristic.

It is known that the Gaussian RBF kernel $\exp(-\|x-y\|^2/\sigma^2)$ and the so-called Laplacian kernel $\exp(-\alpha \sum_{i=1}^{m} |x_i - y_i|)$ $(\alpha > 0)$ are characteristic on $\mathbb{R}^m$ or on a compact subset of $\mathbb{R}^m$ with respect to the Borel $\sigma$-field [2, 15, 28].

The following theorem improves Theorem 7 in [14], and is the theoretical basis of kernel dimension reduction. In the following, let $P_B$ denote the probability on $\mathcal{X}$ induced from $P_X$ by the projection $BB^T : \mathcal{X} \to \mathcal{X}$.

THEOREM 4. *Suppose that the closure of the $\mathcal{H}_{\mathcal{X}}^B$ in $L^2(P_X)$ is included in the closure of $\mathcal{H}_{\mathcal{X}}$ in $L^2(P_X)$ for any $B \in \mathbb{S}_d^m(\mathbb{R})$. Then,*

$$\Sigma_{YY|X}^B \geq \Sigma_{YY|X}, \tag{13}$$

*where the inequality refers to the order of self-adjoint operators. If further $(\mathcal{H}_{\mathcal{X}}, P_X)$ and $(\mathcal{H}_{\mathcal{X}}^B, P_B)$ satisfy (AS) for every $B \in \mathbb{S}_d^m(\mathbb{R})$ and $\mathcal{H}_{\mathcal{Y}}$ is characteristic, the following equivalence holds*

$$\Sigma_{YY|X} = \Sigma_{YY|X}^B \iff Y \perp\!\!\!\perp X | B^T X. \tag{14}$$

PROOF. The first assertion is obvious from Proposition 2. For the second assertion, let $C$ be an $m \times (m-d)$ matrix whose columns span the orthogonal complement to the subspace spanned by the columns of $B$, and let $(U, V) = (B^T X, C^T X)$ for notational simplicity. By taking the expectation of the well-known relation

$$\mathrm{Var}_{Y|U}[g(Y)|U] = E_{V|U}[\mathrm{Var}_{Y|U,V}[g(Y)|U,V]] + \mathrm{Var}_{V|U}[E_{Y|U,V}[g(Y)|U,V]]$$

with respect to $V$, we have

$$E_U[\mathrm{Var}_{Y|U}[g(Y)|U]]$$
$$= E_X[\mathrm{Var}_{Y|X}[g(Y)|X]] + E_U[\mathrm{Var}_{V|U}[E_{Y|U,V}[g(Y)|U,V]]],$$

from which Proposition 3 yields

$$\langle g, (\Sigma_{YY|X}^B - \Sigma_{YY|X})g \rangle_{\mathcal{H}_{\mathcal{Y}}} = E_U[\mathrm{Var}_{V|U}[E_{Y|U,V}[g(Y)|U,V]]].$$

It follows that the right-hand side of the equivalence in (14) holds if and only if $E_{Y|U,V}[g(Y)|U, V]$ does not depend on $V$ almost surely. This is equivalent to

$$E_{Y|X}[g(Y)|X] = E_{Y|U}[g(Y)|U]$$

almost surely. Since $\mathcal{H}_{\mathcal{Y}}$ is characteristic, this means that the conditional probability of $Y$ given $X$ is reduced to that of $Y$ given $U$. □

The assumption (AS) and the notion of characteristic kernel are closely related. In fact, from the following proposition, (AS) is satisfied if a characteristic kernel is used. Thus, if $\mathcal{Y}$ is Euclidean, the choice of Gaussian RBF kernels for $k_d$, $k_\mathcal{X}$ and $k_\mathcal{Y}$ is sufficient to guarantee the equivalence given by (14).

PROPOSITION 5. *Let $(\Omega, \mathcal{B})$ be a measurable space, and $(k, \mathcal{H})$ be a bounded measurable positive definite kernel on $\Omega$ and its RKHS. Then, $k$ is characteristic if and only if $\mathcal{H} + \mathbb{R}$ is dense in $L^2(P)$ for any probability measure $P$ on $(\Omega, \mathcal{B})$.*

PROOF. For the proof of "if" part, suppose $m_P = m_Q$ for $P \neq Q$. Denote the total variation of $P - Q$ by $|P - Q|$. Since $\mathcal{H} + \mathbb{R}$ is dense in $L^2(|P - Q|)$, for arbitrary $\varepsilon > 0$ and $A \in \mathcal{B}$, there exists $f \in \mathcal{H} + \mathbb{R}$ such that $\int |f - I_A| \, d|P - Q| < \varepsilon$, where $I_A$ is the index function of $A$. It follows that $|(E_P[f(X)] - P(A)) - (E_Q[f(X)] - Q(A))| < \varepsilon$. Because $E_P[f(X)] = E_Q[f(X)]$ from $m_P = m_Q$, we have $|P(A) - Q(A)| < \varepsilon$ for any $\varepsilon > 0$, which contradicts $P \neq Q$.

For the opposite direction, suppose $\mathcal{H} + \mathbb{R}$ is not dense in $L^2(P)$. There is nonzero $f \in L^2(P)$ such that $\int f \, dP = 0$ and $\int f \varphi \, dP = 0$ for any $\varphi \in \mathcal{H}$. Let $c = 1/\|f\|_{L^1(P)}$, and define two probability measures $Q_1$ and $Q_2$ by $Q_1(E) = c \int_E |f| \, dP$ and $Q_2(E) = c \int_E (|f| - f) \, dP$ for any measurable set $E$. By $f \neq 0$, we have $Q_1 \neq Q_2$, while $E_{Q_1}[k(\cdot, X)] - E_{Q_2}[k(\cdot, X)] = c \int f(x) k(\cdot, x) \, dP(x) = 0$, which means $k$ is not characteristic. □

2.3. *Kernel dimension reduction procedure.* We now use the characterization given in Theorem 4 to develop an optimization procedure for estimating the central subspace from an empirical sample $\{(X_1, Y_1), \ldots, (X_n, Y_n)\}$. We assume that $\{(X_1, Y_1), \ldots, (X_n, Y_n)\}$ is sampled i.i.d. from $P_{XY}$ and we assume that there exists $B_0 \in \mathbb{S}_d^m(\mathbb{R})$ such that $p_{Y|X}(y|x) = p_{Y|B_0^T X}(y|B_0^T x)$.

We define the *empirical cross-covariance operator* $\widehat{\Sigma}_{YX}^{(n)}$ by evaluating the cross-covariance operator at the empirical distribution $\frac{1}{n}\sum_{i=1}^n \delta_{X_i}\delta_{Y_i}$. When acting on functions $f \in \mathcal{H}_\mathcal{X}$ and $g \in \mathcal{H}_\mathcal{Y}$, the operator $\widehat{\Sigma}_{YX}^{(n)}$ gives the empirical covariance

$$\langle g, \widehat{\Sigma}_{YX}^{(n)} f \rangle_{\mathcal{H}_\mathcal{Y}} = \frac{1}{n}\sum_{i=1}^n g(Y_i) f(X_i) - \left(\frac{1}{n}\sum_{i=1}^n g(Y_i)\right)\left(\frac{1}{n}\sum_{i=1}^n f(X_i)\right).$$

Also, for $B \in \mathbb{S}_d^m(\mathbb{R})$, let $\widehat{\Sigma}_{YY|X}^{B(n)}$ denote the *empirical conditional covariance operator*:

(15) $$\widehat{\Sigma}_{YY|X}^{B(n)} = \widehat{\Sigma}_{YY}^{(n)} - \widehat{\Sigma}_{YX}^{B(n)}(\widehat{\Sigma}_{XX}^{B(n)} + \varepsilon_n I)^{-1} \widehat{\Sigma}_{XY}^{B(n)}.$$

The regularization term $\varepsilon_n I$ ($\varepsilon_n > 0$) is required to enable operator inversion and is thus analogous to Tikhonov regularization [17]. We will see that the regularization term is also needed for consistency.

We now define the KDR estimator $\widehat{B}^{(n)}$ as any minimizer of $\text{Tr}[\widehat{\Sigma}_{YY|X}^{B(n)}]$ on the manifold $\mathbb{S}_d^m(\mathbb{R})$; that is, any matrix in $\mathbb{S}_d^m(\mathbb{R})$ that minimizes

(16) $$\text{Tr}[\widehat{\Sigma}_{YY}^{(n)} - \widehat{\Sigma}_{YX}^{B(n)}(\widehat{\Sigma}_{XX}^{B(n)} + \varepsilon_n I)^{-1}\widehat{\Sigma}_{XY}^{B(n)}].$$

In view of (11), this is equivalent to minimizing

$$\sum_{a=1}^{\infty} \min_{f \in \mathcal{H}_\mathcal{X}^B} \left[ \sum_{i=1}^n \left| \left\{ \xi_a(Y_j) - \frac{1}{n}\sum_{j=1}^n \xi_a(Y_j) \right\} \right.\right.$$
$$\left.\left. - \left\{ f(X_j) - \frac{1}{n}\sum_{j=1}^n f(X_j) \right\} \right|^2 + \varepsilon_n \|f\|_{\mathcal{H}_\mathcal{X}^B}^2 \right]$$

over $B \in \mathbb{S}_d^m(\mathbb{R})$, where $\{\xi_a\}_{a=1}^{\infty}$ is a complete orthonormal system for $\mathcal{H}_\mathcal{Y}$.

The KDR contrast function in (16) can also be expressed in terms of Gram matrices (given a kernel $k$, the *Gram matrix* is the $n \times n$ matrix whose entries are the evaluations of the kernel on all pairs of $n$ data points). Let $\phi_i^B \in \mathcal{H}_\mathcal{X}^B$ and $\psi_i \in \mathcal{H}_\mathcal{Y}$ ($1 \le i \le n$) be functions defined by

$$\phi_i^B = k^B(\cdot, X_i) - \frac{1}{n}\sum_{j=1}^n k^B(\cdot, X_j), \qquad \psi_i = k_\mathcal{Y}(\cdot, Y_i) - \frac{1}{n}\sum_{j=1}^n k_\mathcal{Y}(\cdot, Y_j).$$

Because $\mathcal{R}(\widehat{\Sigma}_{XX}^{B(n)}) = \mathcal{N}(\widehat{\Sigma}_{XX}^{B(n)})^\perp$ and $\mathcal{R}(\widehat{\Sigma}_{YY}^{(n)}) = \mathcal{N}(\widehat{\Sigma}_{YY}^{(n)})^\perp$ are spanned by $(\phi_i^B)_{i=1}^n$ and $(\psi_i)_{i=1}^n$, respectively, the trace of $\widehat{\Sigma}_{YY|X}^{B(n)}$ is equal to that of the matrix representation of $\widehat{\Sigma}_{YY|X}^{B(n)}$ on the linear hull of $(\psi_i)_{i=1}^n$. Note that although the vectors $(\psi_i)_{i=1}^n$ are over-complete, the trace of the matrix representation with respect to these vectors is equal to the trace of the operator.

For $B \in \mathbb{S}_d^m(\mathbb{R})$, the centered Gram matrix $G_X^B$ with respect to the kernel $k^B$ is defined by

$$(G_X^B)_{ij} = \langle \phi_i^B, \phi_j^B \rangle_{\mathcal{H}_\mathcal{X}^B}$$
$$= k_\mathcal{X}^B(X_i, X_j) - \frac{1}{n}\sum_{b=1}^n k_\mathcal{X}^B(X_i, X_b) - \frac{1}{n}\sum_{a=1}^n k_\mathcal{X}^B(X_a, X_j)$$
$$+ \frac{1}{n^2}\sum_{a=1}^n \sum_{b=1}^n k_\mathcal{X}^B(X_a, X_b)$$

and $G_Y$ is defined similarly. By direct calculation, it is easy to obtain

$$\widehat{\Sigma}_{YY|X}^{B(n)} \psi_i = \frac{1}{n}\sum_{j=1}^n \psi_j (G_Y)_{ji} - \frac{1}{n}\sum_{j=1}^n \psi_j (G_X^B(G_X^B + n\varepsilon_n I_n)^{-1} G_Y)_{ji}.$$

It follows that the matrix representation of $\widehat{\Sigma}_{YY|X}^{B(n)}$ with respect to $(\psi_i)_{i=1}^n$ is $\frac{1}{n}\{G_Y - G_X^B(G_X^B + n\varepsilon_n I_n)^{-1}G_Y\}$ and its trace is

$$\operatorname{Tr}[\widehat{\Sigma}_{YY|X}^{B(n)}] = \frac{1}{n}\operatorname{Tr}[G_Y - G_X^B(G_X^B + n\varepsilon_n I_n)^{-1}G_Y]$$
$$= \varepsilon_n \operatorname{Tr}[G_Y(G_X^B + n\varepsilon_n I_n)^{-1}].$$

Omitting the constant factor, the KDR contrast function in (16) thus reduces to

(17) $$\operatorname{Tr}[G_Y(G_X^B + n\varepsilon_n I_n)^{-1}].$$

The KDR method is defined as the optimization of this function over the manifold $\mathbb{S}_d^m(\mathbb{R})$.

Theorem 4 is the population justification of the KDR method. Note that this derivation imposes no strong assumptions either on the conditional probability of $Y$ given $X$, or on the marginal distributions of $X$ and $Y$. In particular, it does not require ellipticity of the marginal distribution of $X$, nor does it require an additive noise model. The response variable $Y$ may be either continuous or discrete. We confirm this general applicability of the KDR method by the numerical results presented in the next section.

Because the contrast function (17) is nonconvex, the minimization requires a nonlinear optimization technique; in our experiments we use the steepest descent method with line search. To alleviate potential problems with local optima, we use a continuation method in which the scale parameter $\sigma$ in Gaussian RBF kernel $\exp(-\|x-y\|/\sigma^2)$ is gradually decreased during the iterative optimization process. In numerical examples shown in the next section, we used a fixed number of iterations, and decreased $\sigma^2$ linearly from $\sigma^2 = 100$ to $\sigma^2 = 10$ for standardized data with standard deviation 5.0. Since the covariance operator approaches the covariance operator induced by a linear kernel as $\sigma \to \infty$, which is solvable as an eigenproblem.

In addition to $\sigma$, there is another tuning parameter $\varepsilon_n$, the regularization coefficient. As both of these tuning parameters have a similar smoothing effect, it is reasonable to fix one of them and select the other; in our experiments we fixed $\varepsilon_n = 0.1$ as an arbitrary choice and varied $\sigma^2$. While there is no theoretical guarantee for this choice, we observe the results are generally stable if the optimization process is successful. There also exist heuristics for choosing kernel parameters in similar RKHS-based dependency analysis; an example is to use the median of pairwise distances of the data for the parameter $\sigma$ in the Gaussian RBF kernel [16]. Currently, however, we are not aware of theoretically justified methods of choosing these parameters; this is an important open problem.

The proposed estimator is shown to be consistent as the sample size goes to infinity. We defer the proof to Section 4.

## 3. Numerical results.

3.1. *Simulation studies.* In this section we compare the performance of the KDR method with that of several well-known dimension reduction methods. Specifically, we compare to SIR, pHd and SAVE on synthetic data sets generated by the regressions in Examples 6.2, 6.3 and 6.4 of [22]. The results are evaluated by computing the Frobenius distance between the projection matrix of the estimated subspace and that of the true subspace; this evaluation measure is invariant under change of basis and is equal to

$$\|B_0 B_0^T - \widehat{B}\widehat{B}^T\|_F,$$

where $B_0$ and $\widehat{B}$ are matrices in the Stiefel manifold $\mathbb{S}_d^m(\mathbb{R})$ representing the true subspace and the estimated subspace, respectively. For the KDR method, a Gaussian RBF kernel $\exp(-\|z_1 - z_2\|^2/c)$ was used, with $c = 2.0$ for regression (A) and regression (C) and $c = 0.5$ for regression (B). The parameter estimate $\widehat{B}$ was updated 100 times by the steepest descent method. The regularization parameter was fixed at $\varepsilon = 0.1$. For SIR and SAVE, we optimized the number of slices for each simulation so as to obtain the best average norm.

Regression (A) is given by

$$\text{(A)} \qquad Y = \frac{X_1}{0.5 + (X_2 + 1.5)^2} + (1 + X_2)^2 + \sigma E,$$

where $X \sim N(0, I_4)$ is a four-dimensional explanatory variable, and $E \sim N(0, 1)$ is independent of $X$. Thus, the central subspace is spanned by the vectors $(1, 0, 0, 0)$ and $(0, 1, 0, 0)$. For the noise level $\sigma$, three different values were used: $\sigma = 0.1, 0.4$ and $0.8$. We used 100 random replications with 100 samples each. Note that the distribution of the explanatory variable $X$ satisfies the ellipticity assumption, as required by the SIR, SAVE and pHd methods.

Table 1 shows the mean and the standard deviation of the Frobenius norm over 100 samples. We see that the KDR method outperforms the other three methods in terms of estimation accuracy. It is also worth noting that in the results presented by Li, Zha and Chiaromonte [22] for their GCR method, the average norm was $0.28, 0.33, 0.45$ for $\sigma = 0.1, 0.4, 0.8$, respectively; again, this is worse than the performance of KDR.

The second regression is given by

$$\text{(B)} \qquad Y = \sin^2(\pi X_2 + 1) + \sigma E,$$

where $X \in \mathbb{R}^4$ is distributed uniformly on the set

$$[0,1]^4 \setminus \{x \in \mathbb{R}^4 | x_i \leq 0.7 \ (i = 1, 2, 3, 4)\},$$

Table 1
*Comparison of KDR and other methods for regression (A)*

|   | KDR | | SIR | | SAVE | | pHd | |
|---|---|---|---|---|---|---|---|---|
| $\sigma$ | NORM | SD | NORM | SD | NORM | SD | NORM | SD |
| 0.1 | 0.11 | 0.07 | 0.55 | 0.28 | 0.77 | 0.35 | 1.04 | 0.34 |
| 0.4 | 0.17 | 0.09 | 0.60 | 0.27 | 0.82 | 0.34 | 1.03 | 0.33 |
| 0.8 | 0.34 | 0.22 | 0.69 | 0.25 | 0.94 | 0.35 | 1.06 | 0.33 |

and $E \sim N(0,1)$ is independent noise. The standard deviation $\sigma$ is fixed at $\sigma = 0.1, 0.2$ and $0.3$. Note that in this example the distribution of $X$ does not satisfy the ellipticity assumption.

Table 2 shows the results of the simulation experiments for this regression. We see that KDR again outperforms the other methods.

The third regression is given by

(C) $$Y = \tfrac{1}{2}(X_1 - a)^2 E,$$

where $X \sim N(0, I_{10})$ is a ten-dimensional variable and $E \sim N(0,1)$ is independent noise. The parameter $a$ is fixed at $a = 0, 0.5$ and $1$. Note that in this example the conditional probability $p(y|x)$ does not obey an additive noise assumption. The mean of $Y$ is zero and the variance is a quadratic function of $X_1$. We generated 100 samples of 500 data.

The results for KDR and the other methods are shown by Table 3, in which we again confirm that the KDR method yields significantly better performance than the other methods. In this case, pHd fails to find the true subspace; this is due to the fact that pHd is incapable of estimating a direction that only appears in the variance [8]. We note also that the results in [22] show that the contour regression methods SCR and GCR yield average norms larger than 1.3.

Although the estimation of variance structure is generally more difficult than that of estimating mean structure, the KDR method nonetheless is effective at finding the central subspace in this case.

Table 2
*Comparison of KDR and other methods for regression (B)*

|   | KDR | | SIR | | SAVE | | pHd | |
|---|---|---|---|---|---|---|---|---|
| $\sigma$ | NORM | SD | NORM | SD | NORM | SD | NORM | SD |
| 0.1 | 0.05 | 0.02 | 0.24 | 0.10 | 0.23 | 0.13 | 0.43 | 0.19 |
| 0.2 | 0.11 | 0.06 | 0.32 | 0.15 | 0.29 | 0.16 | 0.51 | 0.23 |
| 0.3 | 0.13 | 0.07 | 0.41 | 0.19 | 0.41 | 0.21 | 0.63 | 0.29 |

TABLE 3
*Comparison of KDR and other methods for regression (C)*

| $a$ | KDR | | SIR | | SAVE | | pHd | |
|---|---|---|---|---|---|---|---|---|
| | NORM | SD | NORM | SD | NORM | SD | NORM | SD |
| 0.0 | 0.17 | 0.05 | 1.83 | 0.22 | 0.30 | 0.07 | 1.48 | 0.27 |
| 0.5 | 0.17 | 0.04 | 0.58 | 0.19 | 0.35 | 0.08 | 1.52 | 0.28 |
| 1.0 | 0.18 | 0.05 | 0.30 | 0.08 | 0.57 | 0.20 | 1.58 | 0.28 |

3.2. *Applications.* We apply the KDR method to two data sets; one is a binary classification problem and the other is a regression with a continuous response variable. These data sets have been used previously in studies of dimension reduction methods.

The first data set that we studied is *Swiss bank notes* which has been previously studied in the dimension reduction context by Cook and Lee [7], with the data taken from [11]. The problem is that of classifying counterfeit and genuine Swiss bank notes. The data is a sample of 100 counterfeit and 100 genuine notes. There are six continuous explanatory variables that represent aspects of the size of a note: length, height on the left, height on the right, distance of inner frame to the lower border, distance of inner frame to the upper border and length of the diagonal. We standardize each of explanatory variables so that their standard deviation is 5.0.

As we have discussed in the Introduction, many dimension reduction methods (including SIR) are not generally suitable for binary classification problems. Because among inverse regression methods the estimated subspace given by SAVE is necessarily larger than that given by pHd and SIR [7], we compared the KDR method only with SAVE for this data set.

Figure 1 shows two-dimensional plots of the data projected onto the subspaces estimated by the KDR method and by SAVE. The figure shows that the results for KDR appear to be robust with respect to the values of the scale parameter $a$ in the Gaussian RBF kernel. (Note that if $a$ goes to infinity, the result approaches that obtained by a linear kernel, since the linear term in the Taylor expansion of the exponential function is dominant.) In the KDR case, using a Gaussian RBF with scale parameter $a = 10$ and 100 we obtain clear separation of genuine and counterfeit notes. Slightly less separation is obtained for the Gaussian RBF kernel with $a = 10{,}000$, for the linear kernel and for SAVE; in these cases there is an isolated genuine data point that lies close to the class boundary, which is similar to the results using linear discriminant analysis and specification analysis [11]. We see that KDR finds a more effective subspace to separate the two classes than SAVE and the existing analysis. Finally, note that there are two clusters of counterfeit notes in the result of SAVE, while KDR does not show multiple clusters in either class. Although clusters have also been reported in other analyses [11], Section 12, the KDR results suggest that the cluster structure may not be relevant to the classification.

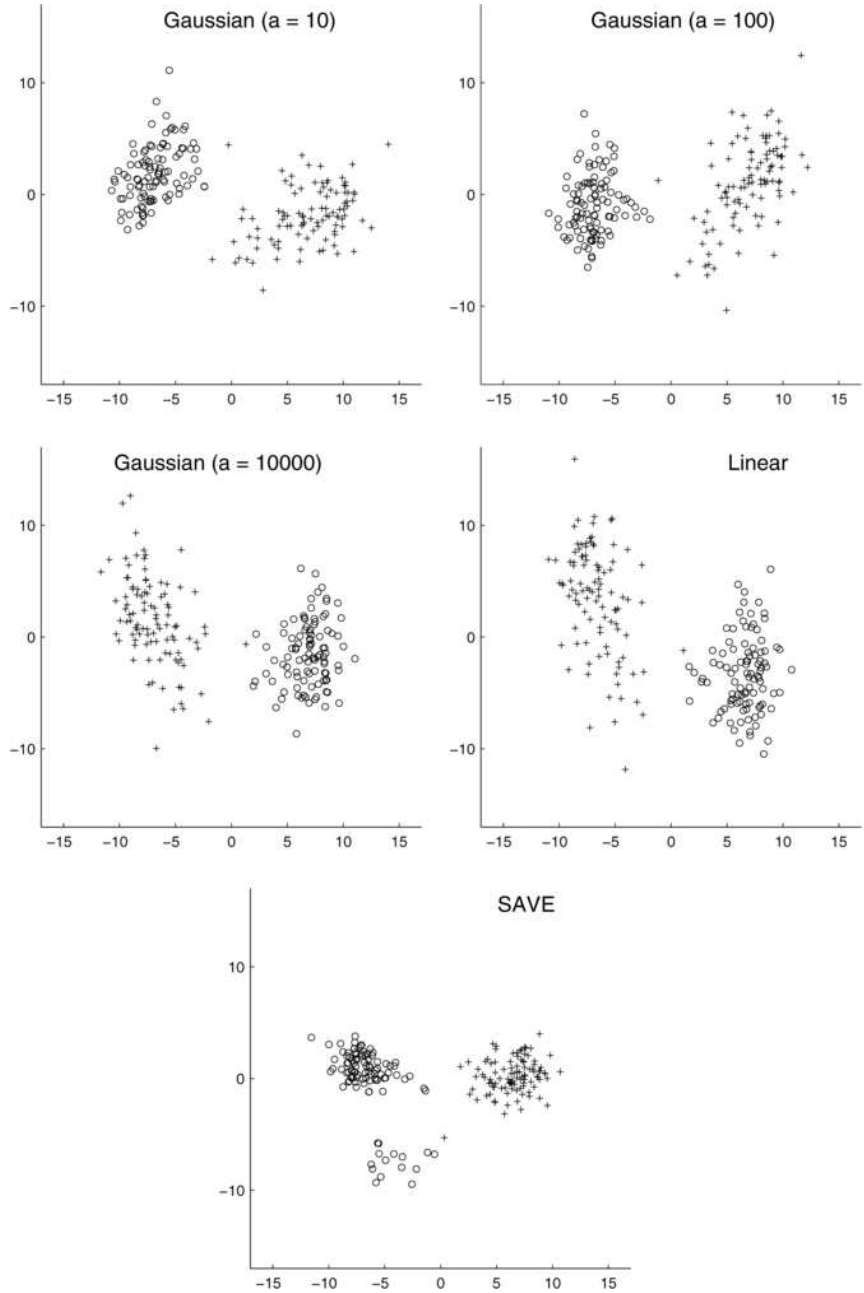

FIG. 1. *Two-dimensional plots of* Swiss bank notes. *The crosses and circles show genuine and counterfeit notes, respectively. For the KDR methods, the Gaussian RBF kernel* $\exp(-\|z_1 - z_2\|^2/a)$ *is used with* $a = 10, 100$ *and* $10{,}000$. *For comparison, the plots given by KDR with a linear kernel and SAVE are shown.*

We also analyzed the *Evaporation* data set, available in the *Arc* package (http://www.stat.umn.edu/arc/software.html). The data set is concerned with the effect on soil evaporation of various air and soil conditions. The number of explanatory variables is 10: maximum daily soil temperature (Maxst), minimum daily soil temperature (Minst), area under the daily soil temperature curve (Avst), maximum daily air temperature (Maxat), minimum daily air temperature (Minat), average daily air temperature (Avat), maximum daily humidity (Maxh), minimum daily humidity (Minh), area under the daily humidity curve (Avh) and total wind speed in miles/hour (Wind). The response variable is daily soil evaporation (Evap). The data were collected daily during 46 days; thus, the number of data points is 46. This data set was studied in the context of contour regression methods for dimension reduction in [22]. We standardize each variable so that the sample variance is equal to 5.0, and use the Gaussian RBF kernel $\exp(-\|z_1 - z_2\|^2/10)$.

Our analysis yielded an estimated two-dimensional subspace which is spanned by the vectors:

$KDR_1$: $-0.25\,MAXST + 0.32\,MINST + 0.00\,AVST + (-0.28)\,MAXAT$
$+ (-0.23)\,MINAT + (-0.44)\,AVAT + 0.39\,MAXH + 0.25\,MINH$
$+ (-0.07)\,AVH + (-0.54)\,WIND.$

$KDR_2$: $0.09\,MAXST + (-0.02)\,MINST + 0.00\,AVST + 0.10\,MAXAT$
$+ (-0.45)\,MINAT + 0.23\,AVAT + 0.21\,MAXH + (-0.41)\,MINH$
$+ (-0.71)\,AVH + (-0.05)\,WIND.$

In the first direction, Wind and Avat have a large factor with the same sign, while both have weak contributions on the second direction. In the second direction, Avh is dominant.

Figure 2 presents the scatter plots representing the response $Y$ plotted with respect to each of the first two directions given by the KDR method. Both of these directions show a clear relation with $Y$. Figure 3 presents the scatter plot of $Y$ versus the two-dimensional subspace found by KDR. The obtained two-dimensional subspace is different from the one given by the existing analysis in [22]; the contour regression method gives a subspace in which the first direction shows a clear monotonic trend, but the second direction suggests a $U$-shaped pattern. In the result of KDR, we do not see a clear folded pattern. Although without further analysis it is difficult to say which result expresses more clearly the statistical dependence, the plots suggest that the KDR method successfully captured the effective directions for regression.

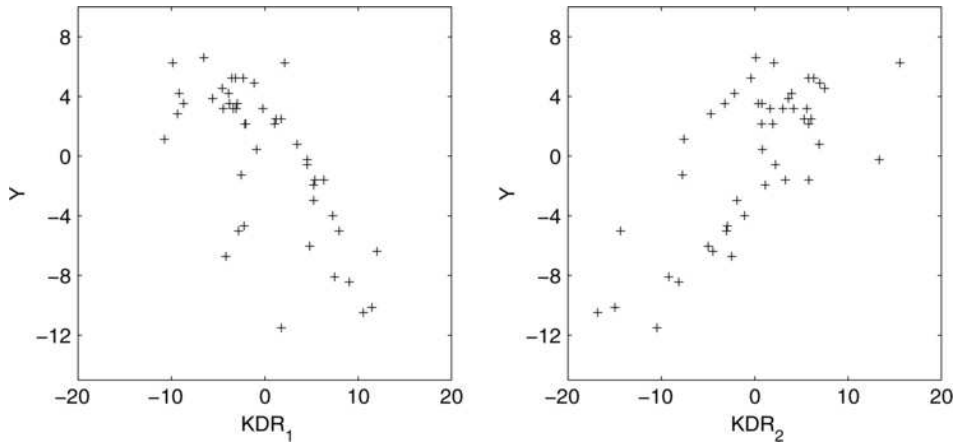

FIG. 2. *Two-dimensional representation of* Evaporation *data for each of the first two directions.*

**4. Consistency of kernel dimension reduction.** In this section we prove that the KDR estimator is consistent. Our proof of consistency requires tools from empirical process theory, suitably elaborated to handle the RKHS setting. We establish convergence of the empirical contrast function to the population contrast function under a condition on the regularization coefficient $\varepsilon_n$, and from this result infer the consistency of $\widehat{B}^{(n)}$.

4.1. *Main result.* We assume hereafter that $\mathcal{Y}$ is a topological space. The Stiefel manifold $\mathbb{S}_d^m(\mathbb{R})$ is assumed to be equipped with a distance $D$ which is compatible with the topology of $\mathbb{S}_d^m(\mathbb{R})$. It is known that geodesics define such a distance (see, e.g., [19], Chapter IV).

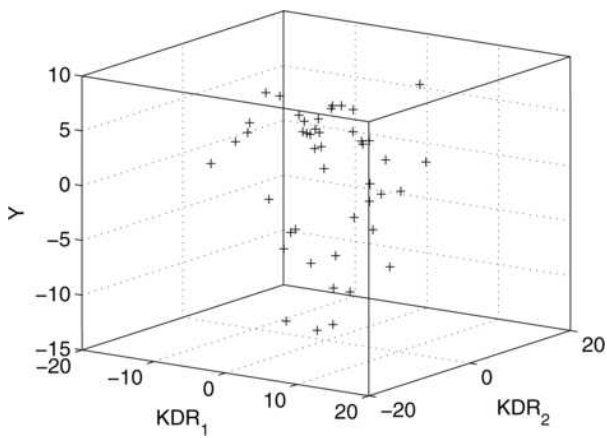

FIG. 3. *Three-dimensional representation of* Evaporation *data.*

The following technical assumptions are needed to guarantee the consistency of kernel dimension reduction:

(A-1) For any bounded continuous function $g$ on $\mathcal{Y}$, the function
$$B \mapsto E_X[E_{Y|B^TX}[g(Y)|B^TX]^2]$$
is continuous on $\mathbb{S}_d^m(\mathbb{R})$.

(A-2) For $B \in \mathbb{S}_d^m(\mathbb{R})$, let $P_B$ be the probability distribution of the random variable $BB^TX$ on $\mathcal{X}$. The Hilbert space $\mathcal{H}_\mathcal{X}^B + \mathbb{R}$ is dense in $L^2(P_B)$ for any $B \in \mathbb{S}_d^m(\mathbb{R})$.

(A-3) There exists a measurable function $\phi: \mathcal{X} \to \mathbb{R}$ such that $E|\phi(X)|^2 < \infty$ and the Lipschitz condition
$$\|k_d(B^Tx, \cdot) - k_d(\tilde{B}^Tx, \cdot)\|_{\mathcal{H}_d} \le \phi(x) D(B, \tilde{B})$$
holds for all $B, \tilde{B} \in \mathbb{S}_d^m(\mathbb{R})$ and $x \in \mathcal{X}$.

THEOREM 6. *Suppose $k_d$ in (9) is continuous and bounded, and suppose the regularization parameter $\varepsilon_n$ in (15) satisfies*

(18) $$\varepsilon_n \to 0, \qquad n^{1/2}\varepsilon_n \to \infty \qquad (n \to \infty).$$

*Define the set of the optimum parameters $\mathbb{B}_d^m$ by*
$$\mathbb{B}_d^m = \underset{B \in \mathbb{S}_d^m(\mathbb{R})}{\arg\min} \operatorname{Tr}[\Sigma_{YY|X}^B].$$

*Under the assumptions* (A-1), (A-2) *and* (A-3), *the set $\mathbb{B}_d^m$ is nonempty, and for an arbitrary open set $U$ in $\mathbb{S}_d^m(\mathbb{R})$ with $\mathbb{B}_d^m \subset U$ we have*
$$\lim_{n \to \infty} \Pr(\widehat{B}^{(n)} \in U) = 1.$$

Note that Theorem 6 holds independently of any requirement that the population contrast function characterizes conditional independence. If the additional conditions of Theorem 4 are satisfied, then the estimator converges in probability to the set of sufficient dimension-reduction subspaces.

The assumptions (A-1) and (A-2) are used to establish the continuity of $\operatorname{Tr}[\Sigma_{YY|X}^B]$ in Lemma 13, and (A-3) is needed to derive the order of uniform convergence of $\widehat{\Sigma}_{YY|X}^{B(n)}$ in Lemma 9.

The assumption (A-1) is satisfied in various cases. Let $f(x) = E_{Y|X}[g(Y)|X = x]$, and assume $f(x)$ is continuous. This assumption holds, for example, if the conditional probability density $p_{Y|X}(y|x)$ is bounded and continuous with respect to $x$. Let $C$ be an element of $\mathbb{S}_{m-d}^m(\mathbb{R})$ such that the subspaces spanned by the column vectors of $B$ and $C$ are orthogonal; that is, the $m \times m$ matrix $(B, C)$ is an orthogonal matrix. Define random variables $U$ and $V$ by $U = B^TX$ and $V = C^TX$. If $X$ has the probability density function $p_X(x)$, the probability density function of $(U, V)$ is given by $p_{U,V}(u,v) = p_X(Bu + Cv)$. Consider the situation in which $u$ is given by $u = B^T \tilde{x}$ for $B \in \mathbb{S}_d^m(\mathbb{R})$ and $\tilde{x} \in \mathcal{X}$, and let $\mathcal{V}_{B,\tilde{x}} = \{v \in \mathbb{R}^{m-d} | BB^T \tilde{x} + Cv \in \mathcal{X}\}$. We have

$$E[g(Y)|B^T X = B^T \tilde{x}] = \frac{\int_{\mathcal{V}_{B,\tilde{x}}} f(BB^T \tilde{x} + Cv) p_X(BB^T \tilde{x} + Cv)\, dv}{\int_{\mathcal{V}_{B,\tilde{x}}} p_X(BB^T \tilde{x} + Cv)\, dv}.$$

If there exists an integrable function $r(v)$ such that $\chi_{\mathcal{V}_{B,\tilde{x}}}(v) p_X(BB^T \tilde{x} + Cv) \leq r(v)$ for all $B \in \mathbb{S}_d^m(\mathbb{R})$ and $\tilde{x} \in \mathcal{X}$, the dominated convergence theorem ensures (A-1). Thus, it is easy to see that a sufficient condition for (A-1) is that $\mathcal{X}$ is bounded, $p_X(x)$ is bounded, and $p_{Y|X}(y|x)$ is bounded and continuous on $x$, which is satisfied by a wide class of distributions.

The assumption (A-2) holds if $\mathcal{X}$ is compact and $k_d + 1$ is a universal kernel on $\mathcal{Z}$. The assumption (A-3) is satisfied by many useful kernels; for example, kernels with the property

$$\left|\frac{\partial^2}{\partial z_a \partial z_b} k_d(z_1, z_2)\right| \leq L\|z_1 - z_2\| \qquad (a, b = 1, 2)$$

for some $L > 0$. In particular Gaussian RBF kernels satisfy this property.

4.2. *Proof of the consistency theorem.* If the following proposition is shown, Theorem 6 follows straightforwardly by standard arguments establishing the consistency of $M$-estimators (see, e.g., [31], Section 5.2).

PROPOSITION 7. *Under the same assumptions as Theorem 6, the functions* $\text{Tr}[\widehat{\Sigma}_{YY|X}^{B(n)}]$ *and* $\text{Tr}[\Sigma_{YY|X}^B]$ *are continuous on* $\mathbb{S}_d^m(\mathbb{R})$*, and*

$$\sup_{B \in \mathbb{S}_d^m(\mathbb{R})} |\text{Tr}[\widehat{\Sigma}_{YY|X}^{B(n)}] - \text{Tr}[\Sigma_{YY|X}^B]| \to 0 \qquad (n \to \infty)$$

*in probability.*

The proof of Proposition 7 is divided into several lemmas. We decompose $\sup_B |\text{Tr}[\Sigma_{YY|X}^B] - \text{Tr}[\widehat{\Sigma}_{YY|X}^{B(n)}]|$ into two parts: $\sup_B |\text{Tr}[\Sigma_{YY|X}^B] - \text{Tr}[\Sigma_{YY} - \Sigma_{YX}^B(\Sigma_{XX}^B + \varepsilon_n I)^{-1} \Sigma_{XY}^B]|$ and $\sup_B |\text{Tr}[\Sigma_{YY} - \Sigma_{YX}^B(\Sigma_{XX}^B + \varepsilon_n I)^{-1} \Sigma_{XY}^B] - \text{Tr}[\widehat{\Sigma}_{YY|X}^{B(n)}]|$. Lemmas 8, 9 and 10 establish the convergence of the second part. The convergence of the first part is shown by Lemmas 11–14; in particular, Lemmas 12 and 13 establish the key result that the trace of the population conditional covariance operator is a continuous function of $B$.

The following lemmas make use of the trace norm and the Hilbert–Schmidt norm of operators. For a discussion of these norms, see [26], Section

VI and [20], Section 30. Recall that the *trace* of a positive operator $A$ on a Hilbert space $\mathcal{H}$ is defined by

$$\mathrm{Tr}[A] = \sum_{i=1}^{\infty} \langle \varphi_i, A\varphi_i \rangle_{\mathcal{H}},$$

where $\{\varphi_i\}_{i=1}^{\infty}$ is a complete orthonormal system (CONS) of $\mathcal{H}$. A bounded operator $T$ on a Hilbert space $\mathcal{H}$ is called *trace class* if $\mathrm{Tr}[(T^*T)^{1/2}]$ is finite. The set of all trace class operators on a Hilbert space is a Banach space with the trace norm $\|T\|_{\mathrm{tr}} = \mathrm{Tr}[(T^*T)^{1/2}]$. For a trace class operator $T$ on $\mathcal{H}$, the series $\sum_{i=1}^{\infty} \langle \varphi_i, T\varphi_i \rangle$ converges absolutely for any CONS $\{\varphi_i\}_{i=1}^{\infty}$, and the limit does not depend on the choice of CONS. The limit is called the *trace* of $T$, and denoted by $\mathrm{Tr}[T]$. It is known that $|\mathrm{Tr}[T]| \leq \|T\|_{\mathrm{tr}}$.

A bounded operator $T : \mathcal{H}_1 \to \mathcal{H}_2$, where $\mathcal{H}_1$ and $\mathcal{H}_2$ are Hilbert spaces, is called *Hilbert–Schmidt* if $\mathrm{Tr}[T^*T] < \infty$, or equivalently, $\sum_{i=1}^{\infty} \|T\varphi_i\|_{\mathcal{H}_2}^2 < \infty$ for a CONS $\{\varphi_i\}_{i=1}^{\infty}$ of $\mathcal{H}_1$. The set of all Hilbert–Schmidt operators from $\mathcal{H}_1$ to $\mathcal{H}_2$ is a Hilbert space with Hilbert–Schmidt inner product

$$\langle T_1, T_2 \rangle_{\mathrm{HS}} = \sum_{i=1}^{\infty} \langle T_1 \varphi_i, T_2 \varphi_i \rangle_{\mathcal{H}_2},$$

where $\{\varphi_i\}_{i=1}^{\infty}$ is a CONS of $\mathcal{H}_1$. Thus, the Hilbert–Schmidt norm $\|T\|_{\mathrm{HS}}$ satisfies $\|T\|_{\mathrm{HS}}^2 = \sum_{i=1}^{\infty} \|T\varphi_i\|_{\mathcal{H}_2}^2$.

Obviously, $\|T\| \leq \|T\|_{\mathrm{HS}} \leq \|T\|_{\mathrm{tr}}$ holds, if $T$ is trace class or Hilbert–Schmidt. Recall also that if $A$ is trace class (Hilbert–Schmidt) and $B$ is bounded, $AB$ and $BA$ are trace class (Hilbert–Schmidt, resp.), for which $\|BA\|_{\mathrm{tr}} \leq \|B\|\|A\|_{\mathrm{tr}}$ and $\|AB\|_{\mathrm{tr}} \leq \|B\|\|A\|_{\mathrm{tr}}$ ($\|AB\|_{\mathrm{HS}} \leq \|A\|\|B\|_{\mathrm{HS}}$ and $\|BA\|_{\mathrm{HS}} \leq \|A\|\|B\|_{\mathrm{HS}}$). If $A : \mathcal{H}_1 \to \mathcal{H}_2$ and $B : \mathcal{H}_2 \to \mathcal{H}_1$ are Hilbert–Schmidt, the product $AB$ is trace-class with $\|AB\|_{\mathrm{tr}} \leq \|A\|_{\mathrm{HS}} \|B\|_{\mathrm{HS}}$.

It is known that cross-covariance operators and covariance operators are Hilbert–Schmidt and trace class, respectively, under the assumption (2) [13, 16]. The Hilbert–Schmidt norm of $\Sigma_{YX}$ is given by

(19) $\quad \|\Sigma_{YX}\|_{\mathrm{HS}}^2 = \|E_{YX}[(k_{\mathcal{X}}(\cdot, X) - m_X)(k_{\mathcal{Y}}(\cdot, Y) - m_Y)]\|_{\mathcal{H}_{\mathcal{X}} \otimes \mathcal{H}_{\mathcal{Y}}}^2,$

where $\mathcal{H}_{\mathcal{X}} \otimes \mathcal{H}_{\mathcal{Y}}$ is the direct product of $\mathcal{H}_{\mathcal{X}}$ and $\mathcal{H}_{\mathcal{Y}}$, and the trace norm of $\Sigma_{XX}$ is

(20) $\quad\quad\quad\quad \mathrm{Tr}[\Sigma_{XX}] = E_X[\|k_{\mathcal{X}}(\cdot, X) - m_X\|_{\mathcal{H}_{\mathcal{X}}}^2].$

LEMMA 8.

$$|\mathrm{Tr}[\widehat{\Sigma}_{YY|X}^{(n)}] - \mathrm{Tr}[\Sigma_{YY} - \Sigma_{YX}(\Sigma_{XX} + \varepsilon_n I)^{-1} \Sigma_{XY}]|$$

$$\leq \frac{1}{\varepsilon_n} \{ (\|\widehat{\Sigma}_{YX}^{(n)}\|_{\mathrm{HS}} + \|\Sigma_{YX}\|_{\mathrm{HS}}) \|\widehat{\Sigma}_{YX}^{(n)} - \Sigma_{YX}\|_{\mathrm{HS}}$$

$$+ \|\Sigma_{YY}\|_{\text{tr}}\|\widehat{\Sigma}_{XX}^{(n)} - \Sigma_{XX}\|\}$$
$$+ |\text{Tr}[\widehat{\Sigma}_{YY}^{(n)} - \Sigma_{YY}]|.$$

PROOF. Noting that the self-adjoint operator $\Sigma_{YX}(\Sigma_{XX} + \varepsilon_n I)^{-1}\Sigma_{XY}$ is trace class from $\Sigma_{YX}(\Sigma_{XX} + \varepsilon_n I)^{-1}\Sigma_{XY} \leq \Sigma_{YY}$, the left-hand side of the assertion is bounded from above by

$$|\text{Tr}[\widehat{\Sigma}_{YY}^{(n)} - \Sigma_{YY}]| + |\text{Tr}[\widehat{\Sigma}_{YX}^{(n)}(\widehat{\Sigma}_{XX}^{(n)} + \varepsilon_n I)^{-1}\widehat{\Sigma}_{XY}^{(n)} - \Sigma_{YX}(\Sigma_{XX} + \varepsilon_n I)^{-1}\Sigma_{XY}]|.$$

The second term is upper-bounded by

$$|\text{Tr}[(\widehat{\Sigma}_{YX}^{(n)} - \Sigma_{YX})(\widehat{\Sigma}_{XX}^{(n)} + \varepsilon_n I)^{-1}\widehat{\Sigma}_{XY}^{(n)}]|$$
$$+ |\text{Tr}[\Sigma_{YX}(\widehat{\Sigma}_{XX}^{(n)} + \varepsilon_n I)^{-1}(\widehat{\Sigma}_{XY}^{(n)} - \Sigma_{XY})]|$$
$$+ |\text{Tr}[\Sigma_{YX}\{(\widehat{\Sigma}_{XX}^{(n)} + \varepsilon_n I)^{-1} - (\Sigma_{XX} + \varepsilon_n I)^{-1}\}\Sigma_{XY}]|$$
$$\leq \|(\widehat{\Sigma}_{YX}^{(n)} - \Sigma_{YX})(\widehat{\Sigma}_{XX}^{(n)} + \varepsilon_n I)^{-1}\widehat{\Sigma}_{XY}^{(n)}\|_{\text{tr}}$$
$$+ \|\Sigma_{YX}(\widehat{\Sigma}_{XX}^{(n)} + \varepsilon_n I)^{-1}(\widehat{\Sigma}_{XY}^{(n)} - \Sigma_{XY})\|_{\text{tr}}$$
$$+ |\text{Tr}[\Sigma_{YX}(\Sigma_{XX} + \varepsilon_n I)^{-1/2}$$
$$\times \{(\Sigma_{XX} + \varepsilon_n I)^{1/2}(\widehat{\Sigma}_{XX}^{(n)} + \varepsilon_n I)^{-1}(\Sigma_{XX} + \varepsilon_n I)^{1/2} - I\}$$
$$\times (\Sigma_{XX} + \varepsilon_n I)^{-1/2}\Sigma_{XY}]|$$
$$\leq \frac{1}{\varepsilon_n}\|\widehat{\Sigma}_{YX}^{(n)} - \Sigma_{YX}\|_{\text{HS}}\|\widehat{\Sigma}_{XY}^{(n)}\|_{\text{HS}} + \frac{1}{\varepsilon_n}\|\Sigma_{YX}\|_{\text{HS}}\|\widehat{\Sigma}_{XY}^{(n)} - \Sigma_{XY}\|_{\text{HS}}$$
$$+ \|(\Sigma_{XX} + \varepsilon_n I)^{1/2}(\widehat{\Sigma}_{XX}^{(n)} + \varepsilon_n I)^{-1}(\Sigma_{XX} + \varepsilon_n I)^{1/2} - I\|$$
$$\times \|(\Sigma_{XX} + \varepsilon_n I)^{-1/2}\Sigma_{XY}\Sigma_{YX}(\Sigma_{XX} + \varepsilon_n I)^{-1/2}\|_{\text{tr}}.$$

In the last line, we use $|\text{Tr}[ABA^*]| \leq \|B\|\|A^*A\|_{\text{tr}}$ for a Hilbert–Schmidt operator $A$ and a bounded operator $B$. This is confirmed easily by the singular decomposition of $A$.

Since the spectrum of $A^*A$ and $AA^*$ are identical, we have

$$\|(\Sigma_{XX} + \varepsilon_n I)^{1/2}(\widehat{\Sigma}_{XX}^{(n)} + \varepsilon_n I)^{-1}(\Sigma_{XX} + \varepsilon_n I)^{1/2} - I\|$$
$$= \|(\widehat{\Sigma}_{XX}^{(n)} + \varepsilon_n I)^{-1/2}(\Sigma_{XX} + \varepsilon_n I)(\widehat{\Sigma}_{XX}^{(n)} + \varepsilon_n I)^{-1/2} - I\|$$
$$\leq \|(\widehat{\Sigma}_{XX}^{(n)} + \varepsilon_n I)^{-1/2}(\Sigma_{XX} - \widehat{\Sigma}_{XX}^{(n)})(\widehat{\Sigma}_{XX}^{(n)} + \varepsilon_n I)^{-1/2}\|$$
$$\leq \frac{1}{\varepsilon_n}\|\widehat{\Sigma}_{XX}^{(n)} - \Sigma_{XX}\|.$$

The bound $\|(\Sigma_{XX} + \varepsilon_n I)^{-1/2}\Sigma_{XX}^{1/2}V_{XY}\| \leq 1$ yields

$$\|(\Sigma_{XX} + \varepsilon_n I)^{-1/2}\Sigma_{XY}\Sigma_{YX}(\Sigma_{XX} + \varepsilon_n I)^{-1/2}\|_{\text{tr}} \leq \|\Sigma_{YY}\|_{\text{tr}},$$

which concludes the proof. □

LEMMA 9. *Under the assumption* (A-3),

$$\sup_{B\in\mathbb{S}_d^m(\mathbb{R})} \|\widehat{\Sigma}_{XX}^{B(n)} - \Sigma_{XX}^B\|_{\mathrm{HS}}, \qquad \sup_{B\in\mathbb{S}_d^m(\mathbb{R})} \|\widehat{\Sigma}_{XY}^{B(n)} - \Sigma_{XY}^B\|_{\mathrm{HS}}$$

*and*

$$\sup_{B\in\mathbb{S}_d^m(\mathbb{R})} |\mathrm{Tr}[\widehat{\Sigma}_{YY}^{B(n)} - \Sigma_{YY}^B]|$$

*are of order* $O_p(1/\sqrt{n})$ *as* $n \to \infty$.

The proof of Lemma 9 is deferred to the Appendix. From Lemmas 8 and 9, the following lemma is obvious.

LEMMA 10. *If the regularization parameter* $(\varepsilon_n)_{n=1}^\infty$ *satisfies* (18), *under the assumption* (A-3) *we have*

$$\sup_{B\in\mathbb{S}_d^m(\mathbb{R})} |\mathrm{Tr}[\widehat{\Sigma}_{YY|X}^{B(n)}] - \mathrm{Tr}[\Sigma_{YY} - \Sigma_{YX}^B(\Sigma_{XX}^B + \varepsilon_n I)^{-1}\Sigma_{XY}^B]| = O_p(\varepsilon_n^{-1} n^{-1/2})$$

*as* $n \to \infty$.

In the next four lemmas, we establish the uniform convergence of $L_\varepsilon$ to $L_0$ ($\varepsilon \downarrow 0$), where $L_\varepsilon(B)$ is a function on $\mathbb{S}_d^m(\mathbb{R})$ defined by

$$L_\varepsilon(B) = \mathrm{Tr}[\Sigma_{YX}^B(\Sigma_{XX}^B + \varepsilon I)^{-1}\Sigma_{XY}^B]$$

for $\varepsilon > 0$ and $L_0(B) = \mathrm{Tr}[\Sigma_{YY}^{1/2} V_{YX}^B V_{XY}^B \Sigma_{YY}^{1/2}]$. We begin by establishing pointwise convergence.

LEMMA 11. *For arbitrary kernels with* (2),

$$\mathrm{Tr}[\Sigma_{YX}(\Sigma_{XX} + \varepsilon I)^{-1}\Sigma_{XY}] \to \mathrm{Tr}[\Sigma_{YY}^{1/2} V_{YX} V_{XY} \Sigma_{YY}^{1/2}] \qquad (\varepsilon \downarrow 0).$$

PROOF. With a CONS $\{\psi_i\}_{i=1}^\infty$ for $\mathcal{H}_\mathcal{Y}$, the difference of the right-hand side and the left-hand side can be written as

$$\sum_{i=1}^\infty \langle \psi_i, \Sigma_{YY}^{1/2} V_{YX}\{I - \Sigma_{XX}^{1/2}(\Sigma_{XX} + \varepsilon I)^{-1}\Sigma_{XX}^{1/2}\} V_{XY} \Sigma_{YY}^{1/2} \psi_i \rangle_{\mathcal{H}_\mathcal{Y}}.$$

Since each summand is positive and upper bounded by $\langle \psi_i, \Sigma_{YY}^{1/2} V_{YX} V_{XY} \Sigma_{YY}^{1/2} \times \psi_i \rangle_{\mathcal{H}_\mathcal{Y}}$, and the sum over $i$ is finite, by the dominated convergence theorem it suffices to show

$$\lim_{\varepsilon \downarrow 0} \langle \psi, \Sigma_{YY}^{1/2} V_{YX}\{I - \Sigma_{XX}^{1/2}(\Sigma_{XX} + \varepsilon I)^{-1}\Sigma_{XX}^{1/2}\} V_{XY} \Sigma_{YY}^{1/2} \psi \rangle_{\mathcal{H}_\mathcal{Y}} = 0$$

for each $\psi \in \mathcal{H}_\mathcal{Y}$.

Fix arbitrary $\psi \in \mathcal{H}_\mathcal{Y}$ and $\delta > 0$. From the fact $\mathcal{R}(V_{XY}) \subset \overline{\mathcal{R}(\Sigma_{XX})}$, there exists $h \in \mathcal{H}_\mathcal{X}$ such that $\|V_{XY}\Sigma_{YY}^{1/2}\psi - \Sigma_{XX}h\|_{\mathcal{H}_\mathcal{X}} < \delta$. Using the fact $I - \Sigma_{XX}^{1/2}(\Sigma_{XX} + \varepsilon_n I)^{-1}\Sigma_{XX}^{1/2} = \varepsilon_n(\Sigma_{XX} + \varepsilon_n I)^{-1}$, we have

$$\|\{I - \Sigma_{XX}^{1/2}(\Sigma_{XX} + \varepsilon I)^{-1}\Sigma_{XX}^{1/2}\}V_{XY}\Sigma_{YY}^{1/2}\psi\|_{\mathcal{H}_\mathcal{X}}$$
$$= \|\varepsilon(\Sigma_{XX} + \varepsilon I)^{-1}\Sigma_{XX}h\|_{\mathcal{H}_\mathcal{X}}$$
$$\quad + \|\varepsilon(\Sigma_{XX} + \varepsilon I)^{-1}(V_{XY}\Sigma_{YY}^{1/2}\psi - \Sigma_{XX}h)\|_{\mathcal{H}_\mathcal{X}}$$
$$\leq \varepsilon\|h\|_{\mathcal{H}_\mathcal{X}} + \delta,$$

which is arbitrarily small if $\varepsilon$ is sufficiently small. This completes the proof. $\square$

LEMMA 12. *Suppose $k_d$ is continuous and bounded. Then, for any $\varepsilon > 0$, the function $L_\varepsilon(B)$ is continuous on $\mathbb{S}_d^m(\mathbb{R})$.*

PROOF. By an argument similar to that in the proof of Lemma 11, it suffices to show the continuity of $B \mapsto \langle \psi, \Sigma_{YX}^B(\Sigma_{XX}^B + \varepsilon I)^{-1}\Sigma_{XY}^B\psi\rangle_{\mathcal{H}_\mathcal{Y}}$ for each $\psi \in \mathcal{H}_\mathcal{Y}$.

Let $J_X^B : \mathcal{H}_X^B \to L^2(P_X)$ and $J_Y : \mathcal{H}_\mathcal{Y} \to L^2(P_Y)$ be inclusions. As seen in Proposition 1, the operators $\Sigma_{YX}^B$ and $\Sigma_{XX}^B$ can be extended to the integral operators $S_{YX}^B$ and $S_{XX}^B$ on $L^2(P_X)$, respectively, so that $J_Y\Sigma_{YX}^B = S_{YX}^B J_X^B$ and $J_X^B\Sigma_{XX}^B = S_{XX}^B J_X^B$. It is not difficult to see also $J_X^B(\Sigma_{XX}^B + \varepsilon I)^{-1} = (S_{XX}^B + \varepsilon I)^{-1} J_X^B$ for $\varepsilon > 0$. These relations yield

$$\langle \psi, \Sigma_{YX}^B(\Sigma_{XX}^B + \varepsilon I)^{-1}\Sigma_{XY}^B\psi\rangle_{\mathcal{H}_\mathcal{Y}}$$
$$= E_{XY}[\psi(Y)((S_{XX}^B + \varepsilon I)^{-1}S_{XY}^B\psi)(X)]$$
$$\quad - E_Y[\psi(Y)]E_X[((S_{XX}^B + \varepsilon I)^{-1}S_{XY}^B\psi)(X)],$$

where $J_Y\psi$ is identified with $\psi$. The assertion is obtained if we prove that the operators $S_{XY}^B$ and $(S_{XX}^B + \varepsilon I)^{-1}$ are continuous with respect to $B$ in operator norm. To see this, let $\tilde{X}$ be identically and independently distributed with $X$. We have

$$\|(S_{XY}^B - S_{XY}^{B_0})\psi\|_{L^2(P_X)}^2$$
$$= E_{\tilde{X}}[\text{Cov}_{YX}[k_\mathcal{X}^B(X, \tilde{X}) - k_\mathcal{X}^{B_0}(X, \tilde{X}), \psi(Y)]^2]$$
$$\leq E_{\tilde{X}}[\text{Var}_X[k_d(B^T X, B^T \tilde{X}) - k_d(B_0^T X, B_0^T \tilde{X})]\text{Var}_Y[\psi(Y)]]$$
$$\leq E_{\tilde{X}}E_X[(k_d(B^T X, B^T \tilde{X}) - k_d(B_0^T X, B_0^T \tilde{X}))^2]\|\psi\|_{L^2(P_Y)}^2,$$

from which the continuity of $B \mapsto S_{XY}^B$ is obtained by the continuity and boundedness of $k_d$. The continuity of $(S_{XX}^B + \varepsilon I)^{-1}$ is shown by $\|(S_{XX}^B +$ $\varepsilon I)^{-1} - (S_{XX}^{B_0} + \varepsilon I)^{-1}\| = \|(S_{XX}^B + \varepsilon I)^{-1}(S_{XX}^{B_0} - S_{XX}^B)(S_{XX}^{B_0} + \varepsilon I)^{-1}\| \leq \frac{1}{\varepsilon^2}\|S_{XX}^{B_0} - S_{XX}^B\|$. $\square$

To establish the continuity of $L_0(B) = \text{Tr}[\Sigma_{YX}^B \Sigma_{XX}^{B}{}^{-1} \Sigma_{XY}^B]$, the argument in the proof of Lemma 12 cannot be applied, because $\Sigma_{XX}^{B}{}^{-1}$ is not bounded in general. The assumptions (A-1) and (A-2) are used for the proof.

LEMMA 13. *Suppose $k_d$ is continuous and bounded. Under the assumptions* (A-1) *and* (A-2), *the function $L_0(B)$ is continuous on $\mathbb{S}_d^m(\mathbb{R})$.*

PROOF. By the same argument as in the proof of Lemma 11, it suffices to establish the continuity of $B \mapsto \langle \psi, \Sigma_{YY|X}^B \psi \rangle$ for $\psi \in \mathcal{H}_\mathcal{Y}$. From Proposition 2, the proof is completed if the continuity of the map

$$B \mapsto \inf_{f \in \mathcal{H}_\mathcal{X}^B} \text{Var}_{XY}[g(Y) - f(X)]$$

is proved for any continuous and bounded function $g$.

Since $f(x)$ depends only on $B^T x$ for any $f \in \mathcal{H}_\mathcal{X}^B$, under the assumption (A-2), we use the same argument as in the proof of Proposition 3 to obtain

$$\inf_{f \in \mathcal{H}_\mathcal{X}^B} \text{Var}_{XY}[g(Y) - f(X)]$$
$$= \inf_{f \in \mathcal{H}_\mathcal{X}^B} \text{Var}_X[E_{Y|BB^T X}[g(Y)|BB^T X] - f(X)]$$
$$\quad + E_X[\text{Var}_{Y|BB^T X}[g(Y)|BB^T X]]$$
$$= E_Y[g(Y)^2] - E_X[E_{Y|B^T X}[g(Y)|B^T X]^2],$$

which is a continuous function of $B \in \mathbb{S}_d^m(\mathbb{R})$ from assumption (A-1). $\square$

LEMMA 14. *Suppose that $k_d$ is continuous and bounded, and that $\varepsilon_n$ converges to zero as $n$ goes to infinity. Under the assumptions* (A-1) *and* (A-2), *we have*

$$\sup_{B \in \mathbb{S}_d^m(\mathbb{R})} \text{Tr}[\Sigma_{YY|X}^B - \{\Sigma_{YY} - \Sigma_{YX}^B(\Sigma_{XX}^B + \varepsilon_n I)^{-1}\Sigma_{XY}^B\}] \to 0 \qquad (n \to \infty).$$

PROOF. From Lemmas 11, 12 and 13, the continuous function $\text{Tr}[\Sigma_{YY} - \Sigma_{YX}(\Sigma_{XX}^B + \varepsilon_n I)^{-1}\Sigma_{XY}^B]$ converges to the continuous function $\text{Tr}[\Sigma_{YY|X}^B]$ for every $B \in \mathbb{S}_d^m(\mathbb{R})$. Because this convergence is monotone and $\mathbb{S}_d^m(\mathbb{R})$ is compact, it is necessarily uniform. $\square$

The proof of Proposition 7 is now easily obtained.

PROOF OF PROPOSITION 7. Lemmas 12 and 13 show the continuity of $\mathrm{Tr}[\widehat{\Sigma}_{YY|X}^{B(n)}]$ and $\mathrm{Tr}[\Sigma_{YY|X}^{B}]$. Lemmas 10 and 14 prove the uniform convergence. $\square$

**5. Conclusions.** This paper has presented KDR, a new method for sufficient dimension reduction in regression. The method is based on a characterization of conditional independence using covariance operators on reproducing Hilbert spaces. This characterization is not restricted to first or second-order conditional moments, but exploits high-order moments in the estimation of the central subspace. The KDR method is widely applicable; in distinction to most of the existing literature on SDR it does not impose strong assumptions on the probability distribution of the covariate vector $X$. It is also applicable to problems in which the response $Y$ is discrete.

We have developed some asymptotic theory for the estimator, resulting in a proof of consistency of the estimator under weak conditions. The proof of consistency reposes on a result establishing the uniform convergence of the empirical process in a Hilbert space. In particular, we have established the rate $O_p(n^{-1/2})$ for uniform convergence, paralleling the results for ordinary real-valued empirical processes.

We have not yet developed distribution theory for the KDR method, and have left open the important problem of inferring the dimensionality of the central subspace. Our proof techniques do not straightforwardly extend to yield the asymptotic distribution of the KDR estimator, and new techniques may be required.

It should be noted, however, that inference of the dimensionality of the central subspace is not necessary for many of the applications of SDR. In particular, SDR is often used in the context of graphical exploration of data, where a data analyst may wish to explore views of varying dimensionality. Also, in high-dimensional prediction problems of the kind studied in statistical machine learning, dimension reduction may be carried out in the context of predictive modeling, in which case cross-validation and related techniques may be used to choose the dimensionality.

Finally, while we have focused our discussion on the central subspace as the object of inference, it is also worth noting that KDR applies even to situations in which a central subspace does not exist. As we have shown, the KDR estimate converges to the subset of projection matrices that satisfy (1); this result holds regardless of the existence of a central subspace. That is, if the intersection of dimension-reduction subspaces is not a dimension-reduction subspace, but if the dimensionality chosen for KDR is chosen to be large enough such that subspaces satisfying (1) exist, then KDR will converge to one of those subspaces.

## APPENDIX: UNIFORM CONVERGENCE OF CROSS-COVARIANCE OPERATORS

In this Appendix we present a proof of Lemma 9. The proof involves the use of random elements in a Hilbert space [3, 30]. Let $\mathcal{H}$ be a Hilbert space equipped with a Borel $\sigma$-field. A *random element* in the Hilbert space $\mathcal{H}$ is a measurable map $F: \Omega \to \mathcal{H}$ from a measurable space $(\Omega, \mathfrak{S})$. If $\mathcal{H}$ is an RKHS on a measurable set $\mathcal{X}$ with a measurable positive definite kernel $k$, a random variable $X$ in $\mathcal{X}$ defines a random element in $\mathcal{H}$ by $k(\cdot, X)$.

A random element $F$ in a Hilbert space $\mathcal{H}$ is said to have *strong order* $p$ ($0 < p < \infty$) if $E\|F\|^p$ is finite. For a random element $F$ of strong order one, the expectation of $F$, which is defined as the element $m_F \in \mathcal{H}$ such that $\langle m_F, g \rangle_\mathcal{H} = E[\langle F, g \rangle_\mathcal{H}]$ for all $g \in \mathcal{H}$, is denoted by $E[F]$. With this notation, the interchange of the expectation and the inner product is justified: $\langle E[F], g \rangle_\mathcal{H} = E[\langle F, g \rangle_\mathcal{H}]$. Note also that for independent random elements $F$ and $G$ of strong order two, the relation

$$E[\langle F, G \rangle_\mathcal{H}] = \langle E[F], E[G] \rangle_\mathcal{H}$$

holds.

Let $(X, Y)$ be a random vector on $\mathcal{X} \times \mathcal{Y}$ with law $P_{XY}$, and let $\mathcal{H}_\mathcal{X}$ and $\mathcal{H}_\mathcal{Y}$ be the RKHS with positive definite kernels $k_\mathcal{X}$ and $k_\mathcal{Y}$, respectively, which satisfy (2). The random element $k_\mathcal{X}(\cdot, X)$ has strong order two, and $E[k(\cdot, X)]$ equals $m_X$, where $m_X$ is given by (4). The random element $k_\mathcal{X}(\cdot, X) k_\mathcal{Y}(\cdot, Y)$ in the direct product $\mathcal{H}_\mathcal{X} \otimes \mathcal{H}_\mathcal{Y}$ has strong order one. Define the zero mean random elements $F = k_\mathcal{X}(\cdot, X) - E[k_\mathcal{X}(\cdot, X)]$ and $G = k_\mathcal{Y}(\cdot, Y) - E[k_\mathcal{Y}(\cdot, Y)]$.

For an i.i.d. sample $(X_1, Y_1), \ldots, (X_n, Y_n)$ on $\mathcal{X} \times \mathcal{Y}$ with law $P_{XY}$, define random elements $F_i = k_\mathcal{X}(\cdot, X_i) - E[k_\mathcal{X}(\cdot, X)]$ and $G_i = k_\mathcal{Y}(\cdot, Y_i) - E[k_\mathcal{Y}(\cdot, Y)]$. Then, $F, F_1, \ldots, F_n$ and $G, G_1, \ldots, G_n$ are zero mean i.i.d. random elements in $\mathcal{H}_\mathcal{X}$ and $\mathcal{H}_\mathcal{Y}$, respectively. In the following, the notation $\mathcal{F} = \mathcal{H}_\mathcal{X} \otimes \mathcal{H}_\mathcal{Y}$ is used for simplicity.

As shown in the proof of Lemma 4 in [13], we have

$$\|\widehat{\Sigma}_{YX}^{(n)} - \Sigma_{YX}\|_{\text{HS}} = \left\| \frac{1}{n} \sum_{i=1}^n \left( F_i - \frac{1}{n} \sum_{j=1}^n F_j \right) \left( G_i - \frac{1}{n} \sum_{j=1}^n G_j \right) - E[FG] \right\|_\mathcal{F},$$

which provides a bound

$$\sup_{B \in \mathbb{S}_d^m(\mathbb{R})} \|\widehat{\Sigma}_{YX}^{B(n)} - \Sigma_{YX}^B\|_{\text{HS}} \leq \sup_{B \in \mathbb{S}_d^m(\mathbb{R})} \left\| \frac{1}{n} \sum_{i=1}^n (F_i^B G_i - E[FG]) \right\|_{\mathcal{F}^B} \quad (21)$$

$$+ \sup_{B \in \mathbb{S}_d^m(\mathbb{R})} \left\| \frac{1}{n} \sum_{j=1}^n F_j^B \right\|_{\mathcal{H}_\mathcal{X}^B} \left\| \frac{1}{n} \sum_{j=1}^n G_j \right\|_{\mathcal{H}_\mathcal{Y}},$$

where $F_i^B$ are defined with the kernel $k^B$, and $\mathcal{F}^B = \mathcal{H}_\mathcal{X}^B \otimes \mathcal{H}_\mathcal{Y}$. Also, (20) implies

$$\mathrm{Tr}[\widehat{\Sigma}_{XX}^{(n)} - \Sigma_{XX}] = \frac{1}{n}\sum_{i=1}^n \left\| F_i - \frac{1}{n}\sum_{j=1}^n F_j \right\|_{\mathcal{H}_\mathcal{X}}^2 - E\|F\|_{\mathcal{H}_\mathcal{X}}^2$$

$$= \frac{1}{n}\sum_{i=1}^n \|F_i\|_{\mathcal{H}_\mathcal{X}}^2 - E\|F\|_{\mathcal{H}_\mathcal{X}}^2 - \left\|\frac{1}{n}\sum_{i=1}^n F_i\right\|_{\mathcal{H}_\mathcal{X}}^2,$$

from which we have

(22)
$$\sup_{B\in\mathbb{S}_d^m(\mathbb{R})} |\mathrm{Tr}[\widehat{\Sigma}_{XX}^{B(n)} - \Sigma_{XX}^B]| \leq \sup_{B\in\mathbb{S}_d^m(\mathbb{R})} \left|\frac{1}{n}\sum_{i=1}^n \|F_i^B\|_{\mathcal{H}_\mathcal{X}^B}^2 - E\|F^B\|_{\mathcal{H}_\mathcal{X}^B}^2\right|$$

$$+ \sup_{B\in\mathbb{S}_d^m(\mathbb{R})} \left\|\frac{1}{n}\sum_{i=1}^n F_i^B\right\|_{\mathcal{H}_\mathcal{X}^B}^2.$$

It follows that Lemma 9 is proved if all the four terms on the right-hand side of (21) and (22) are of order $O_p(1/\sqrt{n})$.

Hereafter, the kernel $k_d$ is assumed to be bounded. We begin by considering the first term on the right-hand side of (22). This is the supremum of a process which consists of real-valued random variables $\|F_i^B\|_{\mathcal{H}_\mathcal{X}^B}^2$. Let $U^B$ be a random element in $\mathcal{H}_d$ defined by

$$U^B = k_d(\cdot, B^T X) - E[k_d(\cdot, B^T X)]$$

and let $C > 0$ be a constant such that $|k_d(z,z)| \leq C^2$ for all $z \in \mathcal{Z}$. From $\|U^B\|_{\mathcal{H}_d} \leq 2C$, we have for $B, \tilde{B} \in \mathbb{S}_d^m(\mathbb{R})$

$$|\|F^B\|_{\mathcal{H}_\mathcal{X}^B}^2 - \|F^{\tilde{B}}\|_{\mathcal{H}_\mathcal{X}^{\tilde{B}}}^2| = |\langle U^B - U^{\tilde{B}}, U^B + U^{\tilde{B}}\rangle_{\mathcal{H}_d}|$$

$$\leq \|U^B - U^{\tilde{B}}\|_{\mathcal{H}_d}\|U^B + U^{\tilde{B}}\|_{\mathcal{H}_d}$$

$$\leq 4C\|U^B - U^{\tilde{B}}\|_{\mathcal{H}_d}.$$

The above inequality, combined with the bound

(23)    $$\|U^B - U^{\tilde{B}}\|_{\mathcal{H}_d} \leq 2\phi(x)D(B,\tilde{B})$$

obtained from assumption (A-3), provides a Lipschitz condition $|\|F^B\|_{\mathcal{H}_\mathcal{X}^B}^2 - \|F^{\tilde{B}}\|_{\mathcal{H}_\mathcal{X}^{\tilde{B}}}^2| \leq 8C\phi(x)D(B,\tilde{B})$, which works as a sufficient condition for the uniform central limit theorem [31], Example 19.7. This yields

$$\sup_{B\in\mathbb{S}_d^m(\mathbb{R})} \left|\frac{1}{n}\sum_{i=1}^n \|F_i^B\|_{\mathcal{H}_\mathcal{X}^B}^2 - E\|F^B\|_{\mathcal{H}_\mathcal{X}^B}^2\right| = O_p(1/\sqrt{n}).$$

Our approach to the other three terms is based on a treatment of empirical processes in a Hilbert space. For $B \in \mathbb{S}_d^m(\mathbb{R})$, let $U_i^B = k_d(\cdot, B^T X_i) - E[k_d(\cdot, B^T X)]$ be a random element in $\mathcal{H}_d$. Then the relation $\langle k^B(\cdot, x), k^B(\cdot, \tilde{x}) \rangle_{\mathcal{H}_\mathcal{X}^B} = k_d(B^T x, B^T \tilde{x}) = \langle k_d(\cdot, B^T x), k_d(\cdot, B^T \tilde{x}) \rangle_{\mathcal{H}_d}$ implies

$$(24) \qquad \left\| \frac{1}{n} \sum_{j=1}^n F_j^B \right\|_{\mathcal{H}_\mathcal{X}^B} = \left\| \frac{1}{n} \sum_{j=1}^n U_j^B \right\|_{\mathcal{H}_d},$$

$$(25) \qquad \left\| \frac{1}{n} \sum_{j=1}^n F_j^B G - E[FG] \right\|_{\mathcal{H}_\mathcal{X}^B \otimes \mathcal{H}_\mathcal{Y}} = \left\| \frac{1}{n} \sum_{j=1}^n U_j^B G - E[U^B G] \right\|_{\mathcal{H}_d \otimes \mathcal{H}_\mathcal{Y}}.$$

Note also that the assumption (A-3) gives

$$(26) \qquad \|U^B G - U^{\tilde{B}} G\|_{\mathcal{H}_d \otimes \mathcal{H}_\mathcal{Y}} \leq 2\sqrt{k_\mathcal{Y}(y,y)} \phi(x) D(B, \tilde{B}).$$

From (23)–(26), the proof of Lemma 9 is completed from the following proposition:

PROPOSITION 15. *Let $(\mathcal{X}, \mathcal{B}_\mathcal{X})$ be a measurable space, let $\Theta$ be a compact metric space with distance $D$, and let $\mathcal{H}$ be a Hilbert space. Suppose that $X, X_1, \ldots, X_n$ are i.i.d. random variables on $\mathcal{X}$, and suppose $F : \mathcal{X} \times \Theta \to \mathcal{H}$ is a Borel measurable map. If $\sup_{\theta \in \Theta} \|F(x; \theta)\|_\mathcal{H} < \infty$ for all $x \in \mathcal{X}$ and there exists a measurable function $\phi : \mathcal{X} \to \mathbb{R}$ such that $E[\phi(X)^2] < \infty$ and*

$$(27) \qquad \|F(x; \theta_1) - F(x; \theta_2)\|_\mathcal{H} \leq \phi(x) D(\theta_1, \theta_2) \qquad (\forall \theta_1, \theta_2 \in \Theta),$$

*then we have*

$$\sup_{\theta \in \Theta} \left\| \frac{1}{\sqrt{n}} \sum_{i=1}^n (F(X_i; \theta) - E[F(X; \theta)]) \right\|_\mathcal{H} = O_p(1) \qquad (n \to \infty).$$

The proof of Proposition 15 is similar to that for a real-valued random process, and is divided into several lemmas.

I.i.d. random variables $\sigma_1, \ldots, \sigma_n$ taking values in $\{+1, -1\}$ with equal probability are called *Rademacher* variables. The following concentration inequality is known for a Rademacher average in a Banach space:

PROPOSITION 16. *Let $a_1, \ldots, a_n$ be elements in a Banach space, and let $\sigma_1, \ldots, \sigma_n$ be Rademacher variables. Then, for every $t > 0$*

$$\Pr\left( \left\| \sum_{i=1}^n \sigma_i a_i \right\| > t \right) \leq 2 \exp\left( -\frac{t^2}{32 \sum_{i=1}^n \|a_i\|^2} \right).$$

PROOF. See [21], Theorem 4.7 and the remark thereafter. □

With Proposition 16, the following exponential inequality is obtained with a slight modification of the standard symmetrization argument for empirical processes.

LEMMA 17. *Let $X, X_1, \ldots, X_n$ and $\mathcal{H}$ be as in Proposition 15, and denote $(X_1, \ldots, X_n)$ by $\mathbf{X}_n$. Let $F: \mathcal{X} \to \mathcal{H}$ be a Borel measurable map with $E\|F(X)\|_{\mathcal{H}}^2 < \infty$. For a positive number $M$ such that $E\|F(X)\|_{\mathcal{H}}^2 < M$, define an event $A_n$ by $\frac{1}{n}\sum_{i=1}^n \|F(X_i)\|^2 \leq M$. Then, for every $t > 0$ and sufficiently large $n$,*

$$\Pr\left(\left\{\mathbf{X}_n \Big| \Big\|\frac{1}{n}\sum_{i=1}^n (F(X_i) - E[F(X)])\Big\|_{\mathcal{H}} > t\right\} \cap A_n\right) \leq 8\exp\left(-\frac{nt^2}{1024M}\right).$$

PROOF. First, note that for any sufficiently large $n$ we have $\Pr(A_n) \geq \frac{3}{4}$ and $\Pr(\|\frac{1}{n}\sum_{i=1}^n(F(X_i) - E[F(X)])\| \leq \frac{t}{2}) \geq \frac{3}{4}$. We consider only such $n$ in the following. Let $\tilde{\mathbf{X}}_n$ be an independent copy of $\mathbf{X}_n$, and let $\tilde{A}_n = \{\tilde{\mathbf{X}}_n | \frac{1}{n}\sum_{i=1}^n \|F(\tilde{X}_i)\|^2 \leq M\}$. The obvious inequality

$$\Pr\left(\left\{\mathbf{X}_n \Big| \Big\|\frac{1}{n}\sum_{i=1}^n (F(X_i) - E[F(X)])\Big\|_{\mathcal{H}} > t\right\} \cap A_n\right)$$
$$\times \Pr\left(\left\{\tilde{\mathbf{X}}_n \Big| \Big\|\frac{1}{n}\sum_{i=1}^n (F(\tilde{X}_i) - E[F(X)])\Big\|_{\mathcal{H}} \leq \frac{t}{2}\right\} \cap \tilde{A}_n\right)$$
$$\leq \Pr\left(\left\{(\mathbf{X}_n, \tilde{\mathbf{X}}_n) \Big| \Big\|\frac{1}{n}\sum_{i=1}^n (F(X_i) - F(\tilde{X}_i))\Big\|_{\mathcal{H}} > \frac{t}{2}\right\} \cap A_n \cap \tilde{A}_n\right)$$

and the fact that $B_n := \{(\mathbf{X}_n, \tilde{\mathbf{X}}_n) | \frac{1}{2n}\sum_{i=1}^n (\|F(X_i)\|^2 + \|F(\tilde{X}_i)\|^2) \leq M\}$ includes $A_n \cap \tilde{A}_n$ gives a symmetrized bound

$$\Pr\left(\left\{\mathbf{X}_n \Big| \Big\|\frac{1}{n}\sum_{i=1}^n (F(X_i) - E[F(X)])\Big\|_{\mathcal{H}} > t\right\} \cap A_n\right)$$
$$\leq 2\Pr\left(\left\{(\mathbf{X}_n, \tilde{\mathbf{X}}_n) \Big| \Big\|\frac{1}{n}\sum_{i=1}^n (F(X_i) - F(\tilde{X}_i))\Big\|_{\mathcal{H}} > \frac{t}{2}\right\} \cap B_n\right).$$

Introducing Rademacher variables $\sigma_1, \ldots, \sigma_n$, the right-hand side is equal to

$$2\Pr\left(\left\{(\mathbf{X}_n, \tilde{\mathbf{X}}_n, \{\sigma_i\}) \Big| \Big\|\frac{1}{n}\sum_{i=1}^n \sigma_i(F(X_i) - F(\tilde{X}_i))\Big\|_{\mathcal{H}} > \frac{t}{2}\right\} \cap B_n\right),$$

which is upper-bounded by

$$4\Pr\left(\left\|\frac{1}{n}\sum_{i=1}^{n}\sigma_i F(X_i)\right\|_{\mathcal{H}} > \frac{t}{4} \text{ and } \frac{1}{2n}\sum_{i=1}^{n}\|F(X_i)\|_{\mathcal{H}}^2 \leq M\right)$$

$$= 4E_{\mathbf{X}_n}\left[\Pr\left(\left\|\frac{1}{n}\sum_{i=1}^{n}\sigma_i F(X_i)\right\|_{\mathcal{H}} > \frac{t}{4}\Big|\mathbf{X}_n\right)1_{\{\mathbf{X}_n \in C_n\}}\right],$$

where $C_n = \{\mathbf{X}_n | \frac{1}{n}\sum_{i=1}^{n}\|F(X_i)\|_{\mathcal{H}}^2 \leq 2M\}$. From Proposition 16, the last line is upper-bounded by $4\exp(-\frac{(nt/4)^2}{32\sum_{i=1}^{n}\|F(X_i)\|^2}) \leq 4\exp(-\frac{nt^2}{1024M})$. □

Let $\Theta$ be a set with semimetric $d$. For any $\delta > 0$, the *covering number* $N(\delta, d, \Theta)$ is the smallest $m \in \mathbb{N}$ for which there exist $m$ points $\theta_1, \ldots, \theta_m$ in $\Theta$ such that $\min_{1 \leq i \leq m} d(\theta, \theta_i) \leq \delta$ holds for any $\theta \in \Theta$. We write $N(\delta)$ for $N(\delta, d, \Theta)$ if there is no confusion. For $\delta > 0$, the *covering integral* $J(\delta)$ for $\Theta$ is defined by

$$J(\delta) = \int_0^\delta (8\log(N(u)^2/u))^{1/2}\,du.$$

The chaining lemma [25], which plays a crucial role in the uniform central limit theorem, is readily extendable to a random process in a Banach space.

LEMMA 18 (Chaining lemma). *Let $\Theta$ be a set with semimetric $d$, and let $\{Z(\theta)|\theta \in \Theta\}$ be a family of random elements in a Banach space. Suppose $\Theta$ has a finite covering integral $J(\delta)$ for $0 < \delta < 1$ and suppose there exists a positive constant $R > 0$ such that for all $\theta, \eta \in \Theta$ and $t > 0$ the inequality*

$$\Pr(\|Z(\theta) - Z(\eta)\| > td(\theta, \eta)) \leq 8\exp\left(-\frac{1}{2R}t^2\right)$$

*holds. Then, there exists a countable subset $\Theta^*$ of $\Theta$ such that for any $0 < \varepsilon < 1$*

$$\Pr\left(\sup_{\theta, \eta \in \Theta^*, d(\theta, \eta) \leq \varepsilon} \|Z(\theta) - Z(\eta)\| > 26RJ(d(\theta, \eta))\right) \leq 2\varepsilon$$

*holds. If $Z(\theta)$ has continuous sample paths, then $\Theta^*$ can be replaced by $\Theta$.*

PROOF. By noting that the proof of the chaining lemma for a real-valued random process does not use any special properties of real numbers but the property of the norm (absolute value) for $Z(\theta)$, the proof applies directly to a process in a Banach space. See [25], Section VII.2. □

PROOF OF PROPOSITION 15. Note that (27) means

$$\left\|\frac{1}{n}\sum_{i=1}^{n}(F(X_i;\theta_1)-F(X_i;\theta_2))\right\|_{\mathcal{H}}^{2} \leq D(\theta_1,\theta_2)^2 \frac{1}{n}\sum_{i=1}^{n}\phi(X_i)^2.$$

Let $M > 0$ be a constant such that $E[\phi(X)^2] < M$, and let $A_n = \{\mathbf{X}_n | \|\frac{1}{n} \times \sum_{i=1}^{n}(F(X_i;\theta_1)-F(X_i;\theta_2))\|_{\mathcal{H}}^2 \leq MD(\theta_1,\theta_2)^2\}$. Since the probability of $A_n$ converges to zero as $n \to \infty$, it suffices to show that there exists $\delta > 0$ such that the probability

$$\mathbb{P}_n = \Pr\left(\mathbf{X}_n | A_n \cap \left\{\sup_{\theta \in \Theta}\left\|\frac{1}{\sqrt{n}}\sum_{i=1}^{n}(F(X_i;\theta)-E[F(X;\theta)])\right\|_{\mathcal{H}} > \delta\right\}\right)$$

satisfies $\limsup_{n \to \infty} \mathbb{P}_n = 0$.

With the notation $\tilde{F}_\theta(x) = F(x;\theta) - E[F(X;\theta)]$, from Lemma 17 we can derive

$$\Pr\left(A_n \cap \left\{\mathbf{X}_n \Big| \left\|\frac{1}{\sqrt{n}}\sum_{i=1}^{n}(\tilde{F}_{\theta_1}(X_i) - \tilde{F}_{\theta_2}(X_i))\right\|_{\mathcal{H}} > t\right\}\right)$$
$$\leq 8\exp\left(-\frac{t^2}{512 \cdot 2MD(\theta_1,\theta_2)^2}\right)$$

for any $t > 0$ and sufficiently large $n$. Because the covering integral $J(\delta)$ with respect to $D$ is finite by the compactness of $\Theta$, and the sample path $\Theta \ni \theta \mapsto \frac{1}{\sqrt{n}}\sum_{i=1}^{n}\tilde{F}_\theta(X_i) \in \mathcal{H}$ is continuous, the chaining lemma implies that for any $0 < \varepsilon < 1$

$$\Pr\left(A_n \cap \left\{\mathbf{X}_n \Big| \sup_{\theta_1,\theta_2 \in \Theta, D(\theta_1,\theta_2) \leq \varepsilon}\left\|\frac{1}{\sqrt{n}}\sum_{i=1}^{n}(\tilde{F}_{\theta_1}(X_i) - \tilde{F}_{\theta_2}(X_i))\right\|_{\mathcal{H}} \right.\right.$$
$$\left.\left. > 26 \cdot 512 M \cdot J(\varepsilon)\right\}\right) \leq 2\varepsilon.$$

Take an arbitrary $\varepsilon \in (0,1)$. We can find a finite number of partitions $\Theta = \bigcup_{a=1}^{\nu(\varepsilon)}\Theta_a$ ($\nu(\varepsilon) \in \mathbb{N}$) so that any two points in each $\Theta_a$ are within the distance $\varepsilon$. Let $\theta_a$ be an arbitrary point in $\Theta_a$. Then the probability $\mathbb{P}_n$ is bounded by

$$\mathbb{P}_n \leq \Pr\left(\max_{1 \leq a \leq \nu(\varepsilon)}\left\|\frac{1}{\sqrt{n}}\sum_{i=1}^{n}\tilde{F}_{\theta_a}(X_i)\right\|_{\mathcal{H}} > \frac{\delta}{2}\right)$$
(28)
$$+ \Pr\left(A_n \cap \left\{\mathbf{X}_n \Big| \sup_{\theta,\eta \in \Theta, D(\theta,\eta) \leq \varepsilon}\left\|\frac{1}{\sqrt{n}}\sum_{i=1}^{n}(\tilde{F}_\theta(X_i) - \tilde{F}_\eta(X_i))\right\|_{\mathcal{H}} > \frac{\delta}{2}\right\}\right).$$

From Chebyshev's inequality the first term is upper-bounded by

$$\nu(\varepsilon)\Pr\left(\left\|\frac{1}{\sqrt{n}}\sum_{i=1}^{n}\tilde{F}_{\theta_a}(X_i)\right\|_{\mathcal{H}} > \frac{\delta}{2}\right) \leq \frac{4\nu(\varepsilon)E\|\tilde{F}_{\theta_a}(X)\|_{\mathcal{H}}^2}{\delta^2}.$$

If we take sufficiently large $\delta$ so that $512MJ(\varepsilon) < \delta/2$ and $\frac{4\nu(\varepsilon)E\|\tilde{F}_{\theta_a}(X)\|_{\mathcal{H}}^2}{\varepsilon} < \delta^2$, the right-hand side of (28) is bounded by $3\varepsilon$, which completes the proof. □

**Acknowledgments.** The authors thank the Editor and anonymous referees for their helpful comments. The authors also thank Dr. Yoichi Nishiyama for his helpful comments on the uniform convergence of empirical processes.

K. Fukumizu
Institute of Statistical Mathematics
4-6-7 Minami-Azabu
Minato-ku, Tokyo 106-8569
Japan
E-mail: fukumizu@ism.ac.jp

F. R. Bach
INRIA—WILLOW Project-Team
Laboratoire d'Informatique
 de l'Ecole Normale Supérieure
CNRS/ENS/INRIA UMR 8548
45, rue d'Ulm 75230 Paris
France
E-mail: francis.bach@mines.org

M. I. Jordan
Department of Statistics
Department of Computer Science
 and Electrical Engineering
University of California
Berkeley, California 94720
USA
E-mail: jordan@stat.berkeley.edu